\documentclass{article}
\usepackage{amsmath,amssymb,amsthm}
\usepackage{graphics}

\newtheorem{theorem}{Theorem}
\newtheorem{lemma}{Lemma}

\newcommand{\conv}{\mathop{\rm conv}\nolimits}
\newcommand{\dist}{\mathop{\rm dist}\nolimits}

\newcommand{\cp}{\mathop{\rm Cap}\nolimits}
\newtheorem{cor}{Corollary}
\newtheorem{prop}{Proposition}
\newtheorem{defn}{Definition}

\title {THE KISSING NUMBER IN FOUR DIMENSIONS}

\author {Oleg R. Musin \thanks{Institute for Mathematical Study of Complex Systems, Moscow State University, Russia omusin@gmail.com}}

\begin{document}
\date{}
\maketitle

\begin{abstract}
The kissing number problem asks for the maximal number $k(n)$ of equal size nonoverlapping
spheres  in $n$-dimensional space  that can touch another sphere of the same size. 
This problem in dimension three was the subject of a famous discussion between Isaac Newton 
and David Gregory in 1694. In three dimensions the problem was finally solved only in 1953 by Sch\"utte and van der Waerden. 

In this paper we present a solution of a long-standing problem   about the kissing number  in four dimensions. Namely, the equality $k(4)=24$ is proved.
The proof is based on a modification of Delsarte's method. 

\end{abstract}

\section {Introduction}

The {\it kissing number} $k(n)$ is the highest number of equal nonoverlapping spheres in ${\bf R}^n$ that can touch another sphere of the same size. In three dimensions the kissing number problem is asking how many white billiard balls can 
{\em kiss} (touch) a black ball. 

The most symmetrical configuration, 12 billiard balls around another, is if the 12 balls are placed at positions corresponding to the vertices of a regular icosahedron concentric with the central ball. However, these 12 outer balls do not kiss each other and may all moved freely. So perhaps if you moved all of them to one side a 13th ball would possibly fit in?       

This problem was the subject of a famous discussion between Isaac Newton 
and David Gregory in 1694. 
It is commonly said that Newton believed the answer was 12 balls, while Gregory 
thought that 13 might be possible. However, Casselman \cite{Cas} found some puzzling features in this story.
 
The Newton-Gregory problem is often called the {\it thirteen spheres problem}.  Hoppe \cite{Hop} thought he had solved the problem in 1874. However, there was a mistake - an analysis of this mistake  was published by Hales \cite{Hales}  in 1994.
Finally, this problem was solved by Sch\"utte and van der Waerden  in 1953 \cite{SvdW2}. A subsequent two-page sketch of a proof was given  by Leech \cite{Lee} in 1956. The thirteen spheres problem continues to be of interest, and several new proofs have been published in the last few years \cite{Hs,Ma,Bor,Ans,Mus2}.

Note that $k(4)\ge 24$.
Indeed, the unit sphere in ${\bf R}^4$ centered at $(0,0,0,0)$ has 24 unit spheres around it, centered at the points $(\pm\sqrt{2},\pm\sqrt{2},0,0)$, with any choice of signs and any ordering of the coordinates. The convex hull of these 24 points yields a famous 4-dimensional regular polytope - the ``24-cell". Its facets are 24 regular octahedra. 

Coxeter proposed upper bounds on $k(n)$ in 1963  \cite{Cox}; for $n=4, 5, 6, 7,$ and 8 these bounds were 26, 48, 85, 146, and 244, respectively. Coxeter's bounds are based on the conjecture that equal size spherical caps on a sphere can be packed no denser than packing where the Delaunay triangulation with vertices at the centers of caps consists of regular simplices. This conjecture has been proved by B\"or\"oczky in 1978 \cite{Bor1}.

The main progress in the kissing number problem in high dimensions was made in the end of 1970s. In 1978: Kabatiansky and Levenshtein have found an asymptotic upper bound $2^{0.401n(1+o(1))}$ for $k(n)$ \cite{Kab}. (Currently known the lower bound is $2^{0.2075n(1+o(1))}$ \cite{Wyn}.) In 1979:
Levenshtein \cite{Lev2}, and independently Odlyzko and  Sloane \cite{OdS} (= 
\cite[Chap.13]{CS}), using Delsarte's method, have proved that
$k(8)=240$, and $k(24)=196560$. This proof is surprisingly short, clean, and technically easier than all proofs in three dimensions. 

However, $n=8, 24$ are the only dimensions in which this method gives a precise result. For other dimensions (for instance, $n=3, 4$) the upper bounds exceed the lower. 
In  \cite{OdS} the Delsarte method was applied in dimensions up to 24 (see \cite[Table 1.5]{CS}). For comparison with the values of Coxeter's bounds on $k(n)$ for $n=4, 5, 6, 7,$ and 8 this method gives 25, 46, 82, 140, and 240, respectively. (For $n=3$  Coxeter's and Delsarte's  methods only gave $k(3)\le 13$ \cite{Cox,OdS}.)

Improvements in the upper bounds on kissing numbers (for $n<24$) were rather weak during next years 
(see \cite[Preface to Third Edition]{CS} for a brief review and references). Arestov and Babenko \cite{AB1} proved that 
the bound $\; k(4)\le25\; $ cannot be improved using Delsarte's method. 
Hsiang \cite{Hs1} claims a proof of $k(4)=24.$ His work has not received yet a positive peer review.

If $M$ unit spheres kiss the unit sphere in ${\bf R}^n$, then the set of kissing points 
is an arrangement on the central sphere such that the (Euclidean) distance between any two points is at least 1. So the kissing number problem can be stated in other way: How many points can be placed on the surface of ${\bf S}^{n-1}$ so that the angular separation between any two points is at least $\pi/3$? 

This leads to an important generalization: a finite subset $X$ of ${\bf S}^{n-1}$ is called a {\it spherical $\psi$-code} if for 
every pair  $(x,y)$ of $X$ the inner product $x\cdot y\le \cos{\psi},$ i.e. the minimal angular separation is at least $\psi.$ Spherical codes have many applications. The main application outside mathematics is in the design of signals for data transmission and storage. There are interesting applications to the numerical evaluation of $n$-dimensional integrals \cite[Chap.3]{CS}.

Delsarte's method (also known in coding theory as Delsarte's linear programming method or Delsarte's scheme) is widely used for finding bounds for codes.
This method is described in 
\cite{CS,Kab} ({see also \cite{PZ} for a beautiful exposition}).

In this paper we present an extension of the Delsarte method
that allowed to prove the bound $k(4)<25$, i.e. $k(4)=24$. This extension yields also a proof for $\; k(3)<13$ \cite{Mus2}.

The first version of these proofs used  numerical solutions of some nonconvex constrained optimization problems \cite{Mus} (see also \cite{PZ}). Now, using geometric approach, we reduced it to relatively simple computations.

The paper is organized as follows: Section 2 shows that the main thorem: $k(4)=24$ easily follows from two lemmas: Lemma A and Lemma B. Section 3 reviews the Delsarte method and gives a proof of Lemma A. Section 4 extends  Delsarte's bounds and reduces the upper bound problem for $\psi$-codes to some optimization problem. Section 5 reduces the dimension of the corresponding optimization problem. Section 6 develops a numerical method for a solution of this optimization problem and gives a proof of Lemma B.

\medskip

\medskip

\noindent{\bf Acknowledgment.} I  wish to thank Eiichi Bannai, 
Dmitry Leshchiner, Sergei Ovchinnikov, Makoto Tagami, G\"unter Ziegler, and especially anonymous referees of this paper for helpful discussions and useful comments.

I am very grateful to Ivan Dynnikov who pointed out a gap in arguments on earlier draft of \cite{Mus}.

\section {The Main Theorem}
Let us introduce the following polynomial of degree nine:\footnote{The polynomial $f_4$ was found by the linear programming method (see details in the Appendix).
This method for $n=4,$ $z=1/2,$ $d=9,$ $N=2000,$ $t_0=0.6058$ gives $E\approx 24.7895.$ For  $f_4$  coefficients were changed to ``better looking" ones with $E\approx 24.8644.$}

$$f_4(t): = \frac{1344}{25}\,t^9 - \frac{2688}{25}\,t^7 + \frac{1764}{25}\,t^5 + \frac{2048}{125}\,t^4 - \frac{1229}{125}\,t^3 - \frac{516}{125}\,t^2 -
 \frac{217}{500}\,t - \frac{2}{125}$$

\medskip

\noindent{\bf Lemma A.} {\em Let $X=\{x_1,\ldots,x_M\}$ be points in the unit sphere ${\bf S}^3$. Then
$$
S(X)=\sum\limits_{i=1}^M\sum\limits_{j=1}^M {f_4(x_i\cdot x_j)}\ge M^2.
$$
  }

\medskip

We give a proof of Lemma A in the next section.

\medskip

\noindent{\bf Lemma B.} {\em Suppose $X=\{x_1,\ldots,x_M\}$ is a subset of ${\bf S}^3$ such that the angular separation between any two distinct points $x_i, x_j$ is at least $\pi/3$. Then
$$
S(X)=\sum\limits_{i=1}^M\sum\limits_{j=1}^M {f_4(x_i\cdot x_j)} < 25M.
$$
  }

\medskip

A proof of Lemma B is given in the end of Section 6.

\medskip

\noindent{\bf Main Theorem.}
$\quad k(4)=24.$
\begin{proof} Let $X$ be a spherical $\pi/3$-code in ${\bf S}^3$ with $M=k(4)$ points.  
Then $X$ satisfies the assumptions in Lemmas A, B.  Therefore, $M^2\le S(X) <25M.$
From this $M<25$ follows, i.e. $M\le 24$. From the other side we have $k(4)\ge 24$, 
showing that $M=k(4)=24.$
\end{proof}

 \section {Delsarte's method}
From here on we will speak of $x\in {\bf S}^{n-1}$ alternatively of points in ${\bf S}^{n-1}$ or of vectors in ${\bf R}^n.$

Let $X = \{x_1, x_2,\ldots, x_M\}$ be any finite subset of the unit sphere ${\bf S}^{n-1}
\subset{\bf R}^n,\\ {\bf S}^{n-1}=\{x: x\in {\bf R}^n,$ $x\cdot x=||x||^2=1\}.$ 
By $\phi_{i,j}=\dist(x_i,x_j)$ we denote the spherical (angular) distance between  $x_i,\,  x_j.$ Clearly,  $ \cos{\phi_{i,j}}=x_i\cdot x_j.$

\medskip

\noindent{\bf 3-A. Schoenberg's theorem.} 
Let $u_1, u_2,\ldots, u_M$ be any real numbers. Then 
$$ ||\sum u_ix_i||^2 = \sum\limits_{i,j} \cos{\phi_{i,j}}u_iu_j \ge 0,$$ or equivalently
the Gram matrix $\Big(\cos{\phi_{i,j}}\Big)$ is  positive semidefinite.

Schoenberg \cite{Scho} extended this property to Gegenbauer polynomials $G_k^{(n)}$.
He proved that {\em the matrix  $\Big(G_k^{(n)}(\cos{\phi_{i,j}})\Big)$  is positive semidefinite for any finite $X\subset {\bf S}^{n-1}$.} 

Schoenberg proved also that the converse holds: {\em if $f(t)$ is a real polynomial and for any finite $X\subset{\bf S}^{n-1}$ the matrix $\big(f(\cos{\phi_{i,j}})\big)$ is positive semidefinite, then $f(t)$ is a linear combination of $G_k^{(n)}(t)$ with nonnegative coefficients.}

\medskip

\noindent{\bf 3-B. The Gegenbauer polynomials.} 
Let us recall definitions of Gegenbauer polynomials. Let polynomials $C_k^{(n)}(t)$ are defined by the expansion
$$(1-2rt+r^2)^{(2-n)/2} = \sum\limits_{k=0}\limits^{\infty}r^kC_k^{(n)}(t).$$ 
Then the polynomials $G_k^{(n)}(t): = C_k^{(n)}(t)/C_k^{(n)}(1)$ are called {\it Gegenbauer} or {\it ultraspherical} polynomials. (So the normalization of $G_k^{(n)}$ is determined by the condition $G_k^{(n)}(1)=1.$) 
Also the Gegenbauer polynomials $G_k^{(n)}$ can be defined by the recurrence formula:
$$G_0^{(n)}=1,\;\; G_1^{(n)}=t,\; \ldots,\; G_k^{(n)}=\frac {(2k+n-4)\,t\,G_{k-1}^{(n)}-(k-1)\,G_{k-2}^{(n)}} {k+n-3}$$

They are orthogonal on the interval $[-1,1]$ with respect to the weight function $\rho(t)=(1-t^2)^{(n-3)/2}$ (see details in \cite{Car,CS,Erd,Scho}). In the case $n=3,\; G_k^{(n)}$ are Legendre polynomials $P_k,$ and $G_k^{(4)}$ are  Chebyshev polynomials of the second kind (but with a different normalization than usual, $U_k(1)=1$),
$$ G_k^{(4)}(t)=U_k(t) = \frac {\sin{((k+1)\phi)}}{(k+1)\sin{\phi}}, \quad t=\cos{\phi}, \quad k=0,1,2,\ldots$$

For instance, $\;\; U_0=1,\;\;\; U_1=t,\;\;\; U_2=(4t^2-1)/3,\;\;\; U_3=2t^3-t,\\ U_4=(16t^4-12t^2+1)/5,\;\ldots,\;  
U_9=(256t^9-512t^7+336t^5-80t^3+5t)/5.$

\medskip

\noindent{\bf 3-C. Delsarte's inequality.}  If a symmetric matrix $A$ is positive semidefinite,  then the sum of all its entries is nonnegative. Schoenberg's theorem implies that
 the matrix $\big(G_k^{(n)}(t_{i,j})\big)$ is positive semidefinite, where $t_{i,j}:=\cos{\phi_{i,j}}, \; $ Then 
$$\sum\limits_{i=1}^M\sum\limits_{j=1}^M {{G_k^{(n)}(t_{i,j})}} \ge 0 \eqno (3.1)$$

\begin{defn} We denote by ${\sf G}_n^+$ the set of continuous functions $f:[-1,1]\to {\bf R}$ representable as series
$$
f(t)=\sum\limits_{k=0}^\infty {c_kG_k^{(n)}(t)}
$$
whose coefficients satisfy the following conditions: 
$$c_0>0,\quad c_k\ge 0 \;  \mbox { for } k=1,2,\ldots, \quad
f(1)=\sum\limits_{k=0}^\infty {c_k}<\infty.$$
\end{defn}

Suppose  $f\in{\sf G}_n^+$ and
let  $$S(X)=S_f(X):=\sum\limits_{i=1}^M\sum\limits_{j=1}^M{f(t_{i,j})}.$$ Using $(3.1),$ we get
$$S(X)=
\sum\limits_{k=0}^\infty c_k\left(\sum\limits_{i=1}^M\sum\limits_{j=1}^M {G_k^{(n)}(t_{i,j})}\right)\ge  
\sum\limits_{i=1}^M\sum\limits_{j=1}^M {c_0G_0^{(n)}(t_{i,j})} =   c_0M^2.$$
Then
$$S(X)\ge c_0M^2. \eqno (3.2)$$

\medskip

\noindent{\bf 3-D. Proof of Lemma A.}
\begin{proof}
The expansion of $f_4$ in terms of $U_k=G_k^{(4)}$ is
$$
f_4 = U_0 + 2\,U_1 + \frac{153}{25}\,U_2 + \frac{871}{250}\,U_3 + 
\frac{128}{25}\,U_4 
+\frac{21}{20}\,U_9
$$
We see that $f_4\in{\sf G}_4^+$ with $c_0=1$. So  Lemma A follows from $(3.2)$.
\end{proof}

\medskip

\noindent{\bf 3-E. Delsarte's bound.} 
Let $X=\{x_1,\ldots,x_M\}\subset {\bf S}^{n-1}$ be a spherical $\psi$-code, i.e. for all $i\neq j,$ $t_{i,j}=\cos{\phi_{i,j}}=x_i\cdot x_j\le z:=\cos{\psi},$ i.e. $t_{i,j}\in [-1,z]$ (but $t_{i,i}=1$). 

Suppose $f\in{\sf G}_n^+$ and
$f(t)\le 0 \; \, \mbox{ for all} \; \; t\in [-1,z],\; $ 
then $\; f(t_{i,j})\le 0\; $ for all $\; i\neq j.\;$ That implies
 $$S_f(X)=Mf(1)+2f(t_{1,2})+\ldots+2f(t_{M-1,M}) \le Mf(1).$$  If we combine this with $(3.2),$ then we get $M \le f(1)/c_0.$ 

{\em Let $A(n,\psi)$ be the maximal size of a $\psi$-code in ${\bf S}^{n-1}$.} Then we have: 
$$A(n,\psi) \le \frac {f(1)}{c_0}\eqno (3.3)$$

The inequality $(3.3)$ play a crucial role in the Delsarte method (see details in  \cite{AB1, AB2, Boyv, CS, Del1, Del2, Kab, Lev2, OdS}). If $z=1/2$ and $c_0=1$,  then $(3.3)$ implies $$k(n)=A(n,\pi/3)\le f(1).$$ Levenshtein \cite{Lev2}, and independently Odlyzko and Sloane \cite{OdS} for $n=8, 24$ have found suitable polynomials $f(t)$: $f(t)\le 0 \; \, \mbox{ for all} \; \; t\in [-1,1/2],\; 
f\in{\sf G}_n^+,\; c_0=1$ with
  $$
f(1)=240\; \, \mbox{for} \; \, n=8;\quad \mbox {and} \quad f(1)=196560 \; \, \mbox{for} \; \,  n=24.
$$ 
Then
$$k(8)\le 240, \quad k(24) \le 196560.$$ For $n=8, \, 24$ the minimal vectors in sphere packings $E_8$ and Leech lattice give these kissing numbers. Thus $k(8)=240,$ and $k(24)=196560.$  

When $n=4,$ a polynomial $f$ of degree 9 with $f(1)\approx 25.5585$ was found in \cite{OdS}. This implies 
$24\le k(4) \le 25.$

\section {An extension of Delsarte's method.}

\noindent{\bf 4-A. An extension of Delsarte's bound.}

 Let $f(t)$ be any real function on the interval $[-1,1]$. Let  for a given $\psi$ $z:=\cos{\psi}$. 
Consider on sphere ${\bf S}^{n-1}$  points $y_0,y_1,\ldots,y_m$  
such that
$$y_i\cdot y_j\le z \; \mbox{ for all }\; i\neq j,\quad f(y_0\cdot y_i)> 0 \; \mbox{ for } \; 1\le i\le m. \eqno (4.1)$$

\begin{defn}  For fixed $y_0\in{\bf S}^{n-1}, m\ge 0, z$, and $f(t)$ let us define the family 
$Q_m(y_0)=Q_m(y_0,n,f)$ of finite sets of points from ${\bf S}^{n-1}$ by the formula
$$
Q_m(y_0):=\left\{
\begin{array}{lcl}
\{y_0\},\qquad \qquad \qquad\qquad\qquad\qquad\qquad\qquad\qquad\qquad\qquad m=0, \\
\{Y =\{y_1,\ldots,y_m\}\subset {\bf S}^{n-1}: \{y_0\}\cup Y \; \mbox{ satisfies } \; (4.1)\},
\quad m\ge 1.
\end{array}
\right.
$$
Denote $\mu=\mu(n,z,f):=\max\{m: Q_m(y_0)\ne\emptyset\}$. \\
For $0\le m\le\mu$ we define the function $H=H_f$ on the family $Q_m(y_0)$:
$$
H(y_0):=f(1)\quad \mbox{for}\quad m=0,
$$
$$
H(y_0;Y)=H(y_0;y_1,\ldots,y_m):=f(1)+f(y_0\cdot y_1)+\ldots+f(y_0\cdot y_m) 
\, \mbox{ for }\, m\ge 1.
$$ 
Let
$$h_m=h_m(n,z,f):=\sup\limits_{Y\in Q_m(y_0)}\{H(y_0;Y)\},\quad
h_{max}:=\max{\{h_0,h_1,\ldots,h_\mu\}}.$$
\end{defn}

\begin{theorem} Suppose $f\in{\sf G}_n^+. \; $ 
Then 
$$ A(n,\psi) \le \frac {h_{max}(n,\cos{\psi},f)}{c_0}=\frac{1}{c_0}\max\{h_0,h_1,\ldots,h_\mu\}.$$
\end{theorem}
\begin{proof} 
 Let $X=\{x_1,\ldots,x_M\}\subset{\bf S}^{n-1}\; $ be a spherical $\psi$-code. Since 
$f\in{\sf G}_n^+,$  $(3.2)$ yields: $S(X)\ge c_0M^2.$ 

 Denote $J(i):=\{j:f(x_i\cdot x_j)> 0, \; j\neq i\}, \; X(i):=\{x_j: j\in J(i)\}.$ 
Then
$$S_i(X):=\sum\limits_{j=1}^M {f(x_i\cdot x_j)}\le f(1)+\sum\limits_{j\in J(i)} {f(x_i\cdot x_j)}=H(x_i;X(i))\le h_{max},$$  
so then $$S(X)=\sum\limits_{i=1}\limits^M S_i(X)\le Mh_{max}.\eqno (4.2)$$

We have $c_0M^2\le S(X) \le Mh_{max},$ i.e. $c_0M\le h_{max}\; $ as required.
\end{proof}

Note that $h_0=f(1).$ If $f(t)\le 0$ for all $t\in [-1,z]$, then 
$\mu(n,z,f)=0,$ i.e. $h_{max}=h_0=f(1).$ Therefore, this theorem yields the Delsarte bound
$M\le f(1)/c_0.$

\medskip

\noindent{\bf 4-B. The class of functions ${\it \Phi}(t_0,z)$.}

The problem of evaluating of $h_{max}$ in general case looks even more complicated than the upper bound problem for spherical $\psi$-codes. It is not clear how to find $\mu$, what is an optimal arrangement for $Y$?
Here we consider this problem only for a very restrictive class of functions ${\it \Phi}(t_0,z)$. For the bound given by Theorem 1 we need $f\in{\sf G}_n^+.$ However, for evaluations of $h_m$ we don't need this assumption. So we are not assume that $f\in{\sf G}_n^+.$

\begin{defn} Let  real numbers $t_0, z$ satisfy $1> t_0>z\ge 0.$ We denote by 
${\it \Phi}(t_0,z)$ the set of  functions $f:[-1,1]\to {\bf R}$ such that
$$f(t)\le 0 \; \mbox{ for } \; t\in [-t_0,z].  $$
\end{defn}

Let $f\in{\it \Phi}(t_0,z)$, and let $Y\in \Omega_m(y_0,n,f)$. Denote 
$$e_0:=-y_0,\quad \theta_0:=\arccos{t_0},\quad \theta_i:=\dist(e_0,y_i) \; \mbox{ for } \; i=1,\ldots,m.$$
(In other words, $e_0$ is the antipodal point to $y_0$.)

It is easy to see that $f(y_0\cdot y_i)>0\; $   only if $\; \theta_i<\theta_0.$
Therefore, \\ {\em $Y$ is a spherical $\psi$-code in the open spherical cap $\cp(e_0,\theta_0)$ of center $e_0$ and radius $\theta_0$ with $\pi/2\ge\psi>\theta_0.$}\\
This assumption is quit restrictive and in particular derives the convexity property for $Y$. 
We are using this property in the next section.

\medskip

\noindent{\bf 4-C. Convexity property.}
A subset of ${\bf S}^{n-1}$ is called {\it spherically convex} if it contains, with every two nonantipodal points, the small arc of the great circle containing them. The closure of a convex set is convex and is the intersection of closed hemispheres (see details in \cite{DGK}).

Let $Y=\{y_1,\ldots,y_m\}\subset \cp(e_0,\theta_0), \; \theta_0<\pi/2$.  Then the convex hull of $Y$ is well defined, and is the intersection of all convex sets containing $Y$.
Denote the convex hull of $Y$ by $\Delta_m=\Delta_m(Y)$.

Recall a definition of a vertex of a convex set:
{\em A point $y\in W$ is called the vertex (extremal point) of a spherically convex closed set $W$, if the set $W\setminus \{y\}$ is spherically convex or, equivalently, there are no points $x,z$ from $W$ for which $y$ is an interior point of the minor arc $\widehat{xz}$ of large radius connecting $x,z$.}

\begin {theorem} Let $Y=\{y_1,\ldots,y_m\}\subset {\bf S}^{n-1}$ be a spherical $\psi$-code.
Suppose $Y\subset\cp(e_0,\theta_0)$, and $0<\theta_0<\psi\le \pi/2.$
Then any $y_k$ is a vertex of $\Delta_m$.
\end {theorem}
\begin{proof} The cases $m=1,2$ are evident. For the case $m=3$ the theorem can be easily proved by contradiction. Indeed, suppose that some point, for instance, $y_2$ is not a vertex of $\Delta_3$. Then, firstly, the set $\Delta_3$  is the arc $\widehat{y_1y_3}$, and, secondly, the point $y_2$ lies on the arc $\widehat{y_1y_3}$. From this it follows that 
$\dist(y_1,y_3)\ge 2\psi$, since $Y$ is a $\psi$-code. From the other hand, according to the triangle inequality, we have
$$2\psi\le\dist(y_1,y_3)\le\dist(e_0,y_1)+\dist(e_0,y_3)<2\theta_0.$$
We obtained the contradiction. It remains to prove the theorem for $m\ge4$.

In this paper we need  only one fact from spherical trigonometry, namely the {\em law of cosines} (or the {\it  cosine theorem}): 
$$\cos{\phi} = \cos{\theta_1}\cos{\theta_2}+\sin{\theta_1}\sin{\theta_2}\cos\varphi,$$ 
where for a spherical triangle $ABC$ the angular lengths of its sides are \\ $\dist(A,B)=\theta_1, \,
\dist(A,C)=\theta_2, \, \dist(B,C)=\phi$, and $\angle{BAC}=\varphi$. 

\medskip

By the assumptions: $$\theta_k=\dist(y_k,e_0)<\theta_0<\psi \; \mbox { for }\; 1\le k\le m; \quad
\phi_{k,j}:=\dist(y_k,y_j)\ge \psi,\; k\ne j.$$
Let us prove that there is no point $y_k$ belonging both to the interior of $\Delta_m$ and relative interior of some facet of dimension $d,\; 1\le d\le\dim{\Delta_m}$. Assume the converse. Then consider  
the great $(n-2)$-sphere $\Omega_k$ such that $y_k\in \Omega_k,$ and $\Omega_k$ is orthogonal to the arc $e_0y_k.$ 
(Note that $\theta_k>0$. Conversely, $y_k=e_0$ and 
$\phi_{k,j}=\theta_j\le\theta_0<\psi.$)

The great sphere $\Omega_k$ divides ${\bf S}^{n-1}$ into two closed hemispheres: $H_1$ and $H_2$. Suppose 
$e_0$ lies in the interior of $H_1$, then at least one $y_j$ belongs $H_2$.
Consider the triangle $e_0y_ky_j$ and denote by $\gamma_{k,j}$ the angle $\angle {e_0y_ky_j}$
in this triangle. The law of cosines yields
$$\cos{\theta_j}=\cos{\theta_k}\cos{\phi_{kj}}+\sin{\theta_k}\sin{\phi_{k,j}}\cos{\gamma_{k,j}}$$
Since $y_j \in H_2,$ we have $\gamma_{k,j}\ge 90^\circ,$ and $\cos{\gamma_{k,j}}\le0$ (Fig. 1). 
From the conditions of Theorem 2 there follow the inequalities
$$
\sin{\theta_k}>0,\quad \sin\phi_{k,j}>0,\quad \cos\theta_k>0,\quad \cos\theta_j>0.$$
Hence, using the cosine theorem we obtain
$$
\cos\theta_j=\cos\theta_k\cos\phi_{k,j}+\sin\theta_k\sin\phi_{k,j}\cos\gamma_{k,j},$$ $$ 
0<\cos\theta_j\le\cos\theta_k\cos\phi_{k,j}.$$
From these inequalities  and $0<\cos\theta_k<1$ there follow that, firstly,
$$
0<\cos\phi_{k,j}\quad \Bigl(\mbox{i.e. } \; \psi\le\phi_{k,j}<\pi/2\Bigr),
$$
and, secondly, the inequalities
$$
\cos\theta_j<\cos\phi_{k,j}\le\cos\psi.
$$
Therefore, $\theta_j>\psi$. This contradiction completes the proof of Theorem 2.
\end{proof}

\begin{center}
\begin{picture}(320,140)(-160,-70)

\put(-110,0){\circle*{3}}
\put(-110,-20){\circle*{3}}
\put(-60,50){\circle*{3}}

\thicklines
\put(-150,0){\line(1,0){70}}

\thinlines
\put(-110,0){\line(1,1){50}}
\put(-110,-20){\line(0,1){20}}

\put(-106,-27){$e_0$}
\put(-119,7){$y_k$}
\put(-75,50){$y_j$}
\put(-145,-20){$H_1$}
\put(-145,30){$H_2$}
\put(-75,-3){$\Omega_k$}

\put(-118,-65){Fig. 1}

\put(80,-65){Fig. 2}

\qbezier(67,54)(90,65)(113,54)
\qbezier (67,54) (55,47) (48,37)
\qbezier (113,54) (125,47) (132,37)

\qbezier (48,37) (38,24) (35,11)
\qbezier (35,11) (31,-4) (30,-20)

\qbezier (132,37) (142,24) (145,11)
\qbezier (145,11) (149,-4) (150,-20)

\put(90,60){\circle*{4}}
\put(92,64){$e_0$}

\put(62,14){\circle*{4}}
\put(121,10){\circle*{4}}

\put(45,-24){\circle*{4}}
\put(135,-24){\circle*{4}}

\thicklines
\qbezier(30,-20)(90,-40)(150,-20)
\thinlines
\qbezier(30,-20)(90,0)(150,-20)

\qbezier(90,60)(67,35)(45,-25)
\qbezier(90,60)(113,35)(135,-25)
\qbezier(62,14)(90,15)(121,10)
\put(30,-37){$\Pi(y_i)$}
\put(130,-37){$\Pi(y_j)$}
\put(51,16){$y_i$}
\put(125,10){$y_j$}
\put(85,3){$\phi_{i,j}$}
\put(85,-39){$\gamma_{i,j}$}
\put(75,30){$\theta_i$}
\put(99,30){$\theta_j$}

\end{picture}
\end{center}

\medskip

\noindent{\bf 4-D. Bounds on $\mu$.}
\begin {theorem} 
Let $Y=\{y_1,\ldots,y_m\}\subset {\bf S}^{n-1}$ be a spherical $\psi$-code.
Suppose $Y\subset\cp(e_0,\theta_0)$, and $\; 0<\psi/2\le\theta_0<\psi\le \pi/2.\; $
Then
$$
m\le A\left(n-1,\arccos{\frac{\cos\psi-\cos^2\theta_0}{\sin^2\theta_0}}\right)
$$ 
\end {theorem}
\begin{proof} 
It is easy to see that the assumption $\; 0<\psi/2\le\theta_0<\psi\le \pi/2\; $  guarantees,  firstly, that the right side of the inequality in Theorem 3 is well defined, secondly, that there is  $Y$ with $m\ge 2$.

If $m\ge 2$, then $y_i\ne e_0.$ Conversely, $\psi\le\dist(y_i,y_j)=\dist(e_0,y_j)=\theta_j<\theta_0$, a contradiction. Therefore,
 the projection $\Pi$ from the pole $e_0$   which sends $x \in {\bf S}^{n-1}$ along its meridian to the equator of the sphere is  defined for all  $y_i$.

Denote  $\gamma_{i,j}:=\dist\left(\Pi(y_i),\Pi(y_j)\right)$ (see Fig. 2).  Then from the law of cosines and the inequality $\cos{\phi_{i,j}}\le z=\cos\psi,$ we get
$$\cos{\gamma_{i,j}}=\frac{\cos{\phi_{i,j}}-\cos{\theta_i}\cos{\theta_j}}{\sin{\theta_i}\sin{\theta_j}}
\le 
\frac{z-\cos{\theta_i}\cos{\theta_j}}{\sin{\theta_i}\sin{\theta_j}}
$$ 
Let
$$R(\alpha,\beta)=\frac{z-\cos{\alpha}\cos{\beta}}{\sin{\alpha}\sin{\beta}}, \; \; \mbox{ then } \; \;
\frac{\partial R(\alpha,\beta)}{\partial\alpha}=\frac{\cos{\beta}-z\cos{\alpha}}{\sin^2{\alpha}\sin{\beta}}.$$ 
We have $\theta_0 < \psi$. Therefore, if $\; 0<\alpha, \beta< \theta_0$, then $ \cos{\beta}> z$. That yields: 
 ${\partial R(\alpha,\beta)}/{\partial\alpha}> 0,$ i.e. $R(\alpha,\beta)$ is a monotone increasing function in $\alpha$. We obtain 
$R(\alpha,\beta)<R(\theta_0,\beta)=R(\beta,\theta_0)<R(\theta_0,\theta_0).$ 

Therefore,
$$\cos{\gamma_{i,j}}\le \frac{z-\cos{\theta_i}\cos{\theta_j}}{\sin{\theta_i}\sin{\theta_j}} <
\frac{z-\cos^2{\theta_0}}{\sin^2{\theta_0}}=\cos\delta.$$
Thus $\Pi(Y)$ is a $\delta$-code on the equator  ${\bf S}^{n-2}$. That yields $m\le A(n-1,\delta)$. 
\end{proof}

\begin{cor} Suppose $f\in {\it \Phi}(t_0,z)$.   If $\; 2t_0^2> z+1,\; $ then
$\; \mu(n,z,f)=1$, otherwise 
$$\mu(n,z,f)\le A\Bigl(n-1,\arccos{\frac{z-t_0^2}{1-t_0^2}}\Bigr),$$ 
\end{cor}
\begin{proof} Let $\cos\psi=z,\; \cos\theta_0=t_0$. Then  $2t_0^2 > z+1$ if and only if 
$\psi>2\theta_0.$ Clearly that in this case the size of any $\psi$-code in the  cap $\cp(e_0,\theta_0)$ is at most 1. Otherwise, $\psi\le 2\theta_0$ and this corollary follows from Theorem 3.
\end{proof}

\begin{cor} Suppose $f\in {\it \Phi}(t_0,z)$.   
Then $$\; \mu(3,z,f)\le 5.$$
\end{cor}

\begin{proof} Note that
$$T=\frac{z-t_0^2}{1-t_0^2}\le\frac{z-z^2}{1-z^2}=\frac{z}{1+z}<\frac{1}{2}. \quad \mbox{Then} \; \; \delta=\arccos{T}>\pi/3.$$
Thus $\; \mu(3,z,f)\le A(2,\delta)\le 2\pi/\delta<6$. 
\end{proof}

\begin{cor} Suppose $f\in {\it \Phi}(t_0,z)$.\\ \\
$(i)\; \, \, $ If $\;t_0> \sqrt{z},\; $ then $\; \mu(4,z,f)\le 4.\\ \\$
$(ii)\;  $ If $\; z=1/2,\; t_0\ge 0.6058,\; $ then 
$\; \mu(4,z,f)\le 6.\\$
\end{cor}

\begin{proof} Denote by $\varphi_k(M)$ the largest angular separation that can be attained in a spherical code on ${\bf S}^{k-1}$ containing $M$ points. In three dimensions the best codes and the values $\varphi_3(M)$ presently known for $M\le 12$   and $M=24$ (see \cite{Dan,FeT,SvdW1}).
Sch\"utte and van der Waerden \cite{SvdW1} proved that  $$\varphi_3(5)=\varphi_3(6)=90^\circ, \quad
\cos{\varphi_3(7)}=\cot{40^\circ}\cot{80^\circ}, \quad\varphi_3(7)\approx 77.86954^\circ.$$ 
$(i)$  Since $z-t_0^2<0$, Corollary 1 yields: $\mu(4,z,f)\le A(3,\delta)$, where $\delta>90^\circ.$ We have $\delta>\varphi_3(5).\; $ Thus $\; \mu<5.$

\noindent $(ii)$ Note that for $\; t_0\ge 0.6058,$ $$\arccos{\frac{1/2-t_0^2}{1-t_0^2}}>77.87^\circ.$$
So Corollary 1 implies $\mu(4,1/2,f)\le A(3,77.87^\circ).$
Since $\; 77.87^\circ>\varphi_3(7),\; $ we have $\; A(3,77.87^\circ)<7,\; $ i.e. $\; \mu\le 6.$
\end{proof}

\medskip

\noindent{\bf 4-E. Optimization problem.} Let  
$$t_0:=\cos\theta_0,\quad z:=\cos\psi, \quad
\cos\delta:=\frac{z-t_0^2}{1-t_0^2}, \quad \mu^*:=A(n-1,\delta).$$
For given $n, \psi, \theta_0, f\in {\it \Phi}(t_0,z), e_0\in {\bf S}^{n-1},$ and $m\le \mu^*$, 
the value $h_{m}(n,z,f)$ is the solution of the following optimization problem on ${\bf S}^{n-1}$: 
$$
\mbox{\em maximize } \; f(1)+f(-e_0\cdot y_1)+\ldots+f(-e_0\cdot y_m) $$ 
{\em subject to the constraints}
$$
y_i\in {\bf S}^{n-1}, \; i=1,\ldots,m, \quad \dist(e_0,y_i)\le\theta_0, \quad \dist(y_i,y_j)\ge \psi, \; i\ne j.
$$

The dimension of this problem is $(n-1)m\le(n-1)\mu^*.$
If $\mu^*$ is small enough, then  for small $n$ it gets relatively small -
dimensional optimization problems for computation of values $h_m$. 
If additionally $f(t)$ is a monotone decreasing function on $[-1,-t_0]$, then in some cases 
this problem can be reduced to $(n-1)$ - dimensional optimization problem of a type that can be  treated numerically.

\section{Optimal and irreducible sets}
{\bf 5-A. The monotonicity assumption and optimal sets.}
\begin{defn}  We denote by 
${\it \Phi}^*(z)$ the set of all functions  
$f\in \bigcup\limits_{\tau_0>z}{\it \Phi}(\tau_0,z)$  
such that 
$\; f(t) \; \mbox{ is a monotone decreasing function on the interval }\; [-1,-\tau_0],$
and $f(-1)>0>f(-\tau_0).$

For any $f\in {\it \Phi}^*(z)$, denote $t_0=t_0(f):=\sup\{t\in [\tau_0,1]: f(-t)<0\}$.
\end{defn}

Clearly, if $f\in {\it \Phi}^*(z)$, then  
$f\in {\it \Phi}(t_0,z)$, i.e. $f(t)\le 0$ for $t\in [-1,-t_0]$. Moreover, if $f(t)$    
is a continious function on  $[-1,-z]$, then $f(-t_0)=0$.

Consider a spherical $\psi$-code  $Y=\{y_1,\ldots,y_m\}\subset\cp(e_0,\theta_0)\subset {\bf S}^{n-1}$.
Then we have the constraint: $\phi_{i,j}:=\dist(y_i,y_j)\ge \psi$ for all $i\ne j.$ {\em Denote by
$\Gamma_\psi(Y)$ the graph with the set of vertices $Y$ and the set of edges $y_iy_j$ with  $\phi_{i,j}=\psi.$}

\begin{defn} Let $f\in {\it \Phi}^*(z), \; \psi=\arccos(z), \; \theta_0=\arccos(t_0).$
 We say that a spherical $\psi$-code  $Y=\{y_1,\ldots,y_m\}\subset\cp(e_0,\theta_0)\subset {\bf S}^{n-1}$ is  optimal for $f$ if $\; H_f(-e_0;Y)=h_m(n,z,f).\; $ 

If optimal $Y$ is not unique up to isometry, then we call $Y$ as optimal if the graph $\Gamma_\psi(Y)$ has the maximal number of  edges.
\end{defn}

Let  $\theta_k:=\dist(y_k,e_0)$. Then $H(-e_0;Y)$ can be represented in the form:
$$F_f(\theta_1,\ldots,\theta_m):=H_f(-e_0;Y)=f(1)+f(-\cos{\theta_1})+\ldots+f(-\cos{\theta_m}). $$

Let us call $F(\theta_1,\ldots,\theta_m)=F_f(\theta_1,\ldots,\theta_m)$ {\em the efficient function}. Clearly, if $f\in {\it \Phi}^*(z),$ then the efficient function is {\em a monotone decreasing function in the interval $[0,\theta_0]$ for any variable $\theta_k$}.

\medskip

\noindent{\bf 5-B. Irreducible sets.}
\begin{defn} Let $0<\theta_0<\psi\le\pi/2$.  We say that a spherical $\psi$-code  
$Y=\{y_1,\ldots,y_m\}\subset\cp(e_0,\theta_0)\subset {\bf S}^{n-1}$ is 
irreducible (or jammed) if any $y_k$ can not be shifted towards $e_0$ (i.e. this shift decreases $\theta_k$) such that $Y'$, which is obtained after this shifting, is also a $\psi$-code. 

As above, in the case when irreducible $Y$ is not defined uniquely up to isometry by $\theta_i$, we say that $Y$ is irreducible  if the graph $\Gamma_\psi(Y)$ has the maximal number of  edges.
\end{defn}

\begin{prop} Let $f\in {\it \Phi}^*(z)$. Suppose $Y\subset\cp(e_0,\theta_0)\subset {\bf S}^{n-1}$ is optimal for $f$.  Then $Y$ is irreducible.
\end{prop}
\begin{proof}
The efficient function $F(\theta_1,\ldots,\theta_m)$  increases whenever $\theta_k$ decreases. From this follows that  $y_k$ can not be shifted towards $e_0.$ In the converse case, 
$H(-e_0;Y)=F(\theta_1,\ldots,\theta_m)$ increases whenever $y_k$ tends to $e_0$. It contradicts  the optimality of the initial set $Y$.
\end{proof}

\begin{lemma} If $Y=\{y_1,\ldots,y_m\}$ is irreducible, then\\
$(i)\; e_0\in \Delta_m=$convex hull of ${Y};\\$
$ (ii)$ If $m>1$, then $\deg{y_i}>0$ for all $y_i\in Y$, where by $\deg{y_i}$ denoted the degree of the vertex $y_i$ in the graph  $\Gamma_\psi(Y)$. 
\end{lemma}
\begin{proof}

\noindent $(i)$ Otherwise whole $Y$ can be shifted towards $e_0.$

\noindent $(ii)$ Clearly, if $\phi_{i,j}>\psi$ for all $j\ne i$, then $y_i$ can be shifted towards $e_0$. 
\end{proof}

For $m=1$ from this follows that $e_0=y_1$, i.e. $h_1=\sup\{F(\theta_1)\}=F(0).$ Thus
$$h_1=f(1)+f(-1). \eqno (5.1)$$

For $m=2$, Lemma 1  implies that $\dist{(x_1,x_2)}=\psi$, i.e.
$$\Delta_2=y_1y_2 \; \mbox{\em is an arc of length } \psi.\eqno (5.2)$$

Consider $\Delta_m\subset {\bf S}^{n-1}$ of dimension $k, \; \dim{\Delta_m}=k$. Since $\Delta_m$ is a convex set, there exists the great $k$-dimensional sphere ${\bf S}^k$ in ${\bf S}^{n-1}$ containing 
$\Delta_m.$ 

Note that if $\dim{\Delta_m}=1$, then $m=2.$ Indeed, since $\dim{\Delta_m}=1$, it follows that $Y$ belongs to the great circle ${\bf S}^1$. It is clear that in this case $m=2$. (For instance, $m>2$ contradicts Theorem 2 for $n=2$.)

\medskip

To prove our main results in this section for $n=3,4$ we need the following fact. (For $n=3$, when $\Delta$ is an arc, a proof of this claim  is trivial.) 
\begin{lemma} Consider in ${\bf S}^{n-1}$ an arc  $\omega$ and a regular simplex $\Delta$, both are with edge lengths $\psi,\; \psi\le \pi/2$. Suppose the intersection of $\omega$ and $\Delta$  is not empty. Then at least one of the distances between vertices of $\omega$ and $\Delta$  is less than $\psi$. 
\end{lemma}
\begin{proof} We have $\omega=u_1u_2,\; \Delta=v_1v_2\ldots v_k,\; \dist(u_1,u_2)=\dist(v_i,v_j)=\psi.$

Assume the converse. Then $\dist(u_i,v_j)\ge\psi$ for all $i, j.$ By $U$ denote the union of the spherical caps of centers $v_i,\; i=1,\ldots,k,\; $ and radius $\psi.$ Let $B$ be the boundary of $U.$
Note that $u_1$ and $u_2$ don't lie inside $U.$ If $\{u_1', u_2'\}=\omega\bigcap B$, then
$\psi=\dist(u_1,u_2)\ge\dist(u_1',u_2')$, and $\omega'\bigcap\Delta\ne \emptyset,$  where $\omega'=u_1'u_2'.$

We have the following optimization problem: to find an arc $w_1w_2$ of minimal length subject to the constraints 
$w_1, w_2 \in B$, and $w_1w_2\bigcap\Delta\neq\emptyset$? It is not hard to prove that $\dist(w_1,w_2)$ attains its minimum when $w_1$ and $w_2$ are at the distance of $\psi$ from all $v_i$, i.e. 
$w_1v_1\ldots v_k$ and $w_2v_1\ldots v_k$ are regular simplices with the common facet $\Delta$.
Using this, it can be shown by direct calculation that
$$\cos{\alpha}=\frac{2kz^2-(k-1)z-1}{1+(k-1)z},\quad \alpha=\min{\dist(w_1,w_2)},\; z=\cos{\psi}\eqno (5.3)$$

We have $\alpha\le\psi$. From $(5.3)$ follows that $\cos{\alpha}\ge z$ if and only if $z\ge 1$ or $(k+1)z+1\le 0$. It contradicts the assumption $0\le z<1.$
\end{proof}

\medskip

\noindent{\bf 5-C. Irreducible sets in ${\bf S}^2$.}
Now we consider irreducible sets for $n=3$. In this case $\dim{\Delta_m\le 2}.$

\begin{theorem} Suppose $Y$ is irreducible and $\; \dim(\Delta_m)=2.\\$ Then $3\le m\le 5$, and $\Delta_m$ is a spherical regular triangle, rhomb, or equilateral pentagon with edge lengths $\psi.$
\end{theorem}
\begin{proof} From Corollary 2 follows that $m\le 5.$ On the other hand, $m>2.$ Then
$m=3,\, 4,\, 5.$
Theorem 2 implies that $\Delta_m$ is a convex polygon with vertices $y_1,\ldots,y_m$. 
From Lemma 1 it follows that $e_0\in \Delta_m$, and $\deg{y_i}\geqslant 1.$

First let us prove that if $\deg{y_i}\ge 2$ for all $i$, then  $\Delta_m$ is an equilateral $m$-gon with edge lengths $\psi.$ Indeed, it is clear for $m=3.$

 Lemma 2 implies that two diagonals of $\Delta_m$ of lengths $\psi$ do not intersect each other.
That yields the proof for $m=4.$
When $m=5$, it remains to consider the case where $\Delta_5$ consists of two regular non overlapping triangles with a common vertex (Fig. 3). This case contradicts the convexity of $\Delta_5$. Indeed, since the angular sum in spherical triangle is strictly greater than $180^\circ$ and a larger side of spherical triangle subtends opposite large angle, we have
$\angle{y_iy_1y_j}>60^\circ$. Then 
$$180^\circ\ge \angle{y_2y_1y_5}=\angle{y_2y_1y_3}+\angle{y_3y_1y_4}+\angle{y_4y_1y_5}>180^\circ$$ - a contradiction.

Now we prove that $\deg{y_i}\ge 2.$
Suppose  $\deg{y_1}=1,$ i.e.  $\phi_{1,2}=\psi, \; \, \phi_{1,i}>\psi \, $ for $\, i=3,\ldots,m.$ (Recall that $\phi_{i,j}=\dist(y_i,y_j)$.)
If $e_0\notin y_1y_2$, then after sufficiently small turn of $y_1$ round $y_2$ to $e_0$ (Fig. 4)  the distance $\theta_1$ decreases - a contradiction.  (This turn will be considered in Lemma 3 with more details.)

\begin{center}
\begin{picture}(320,140)(-160,-70)
\put(-112,44){$y_3$}
\put(-90,-30){$y_1$}
\put(-56,44){$y_4$}
\put(-150,6){$y_2$}
\put(-18,6){$y_5$}

\put(-80,-20){\circle*{4}}
\put(-60,40){\circle*{4}}
\put(-140,0){\circle*{4}}
\put(-100,40){\circle*{4}}
\put(-20,0){\circle*{4}}

\put(-80,-20){\line(-3,1){60}}
\put(-80,-20){\line(-1,3){20}}
\put(-80,-20){\line(1,3){20}}
\put(-80,-20){\line(3,1){60}}

\put(-60,40){\line(1,-1){40}}
\put(-140,0){\line(1,1){40}}

\put(-95,-65){Fig. 3}

\put(80,-65){Fig. 4}

\put(40,30){\circle*{4}}
\put(90,-20){\circle*{6}}
\put(140,30){\circle*{4}}
\put(80,20){\circle*{4}}

\thicklines
\put(90,-20){\line(-1,1){50}}
\thinlines
\put(90,-20){\line(1,1){50}}
\put(40,30){\vector(1,1){13}}

\put(98,-25){$y_2$}
\put(86,18){$e_0$}
\put(25,28){$y_1$}
\put(147,28){$y_3$}

\end{picture}
\end{center}

It remains to consider the case: $e_0\in y_1y_2.$
If $\phi_{i,j}=\psi$ where $i>2$ or $j>2$, then $e_0\notin y_iy_j$. Indeed, in the converse case, we have two intersecting diagonals of lengths $\psi.$ Therefore, 
 $\deg{y_i}\ge 2$ for $2<i\le m$. For $m=3, 4$ it implies the proof. For $m=5$ there is the case where $Q_3=y_3y_4y_5$ is a regular triangle of side length $\psi.$ Note that $y_1y_2$ cannot intersect $Q_3$ (otherwise we again have intersecting diagonals of lengths $\psi$), then $y_1y_2$ is a side of $\Delta_5$. In this case, as above, after sufficiently small turn of $Q_3$ round $y_2$ to $e_0$ the distance $\theta_i, \; i=3,4,5,$ decreases - a contradiction. 
\end{proof}

\medskip

\noindent{\bf 5-D. Rotations and irreducible sets in $n$ dimensions.}
Now we extend these results to $n$ dimensions.\footnote{ In the first version of this paper for $m \ge n$ has been claimed that any vertex of $\Gamma_\psi(Y)$ has degree at least $n-1$. 
However,  E. Bannai, M. Tagami, and referees of this paper found some gaps in our exposition. Most of them are related to ``degenerated" configurations. 
In this paper we need only the case $n=4,  m<6$. For this case  Bannai and  Tagami  verified each step of our proof, considered all  ``degenerated" configurations, and finally gave clean and detailed proof (see E. Bannai and M. Tagami: On optimal sets in Musin's 
paper ``The kissing number in four dimensions" in the 
Proceedings of the COE Workshop on Sphere Packings, November 
1-5, 2004,  in Fukuoka Japan).   
Now this claim for all $n$ can be considered only as conjecture. In {\bf 5-D} we prove the claim when $\{y_i\}$ are in ``general position". 
I wish to thank Eiichi Bannai, Makoto Tagami, and anonymous  referees for helpful and useful comments.} 

Let us consider a rotation $R(\varphi,\Omega)$ on ${\bf S}^{n-1}$ about an $(n-3)$ - dimensional great
sphere $\Omega$ in ${\bf S}^{n-1}.$ Without loss of generality, we may assume that
$$\Omega=\{\vec u=(u_1,\ldots,u_n)\in {\bf R}^n: u_1=u_2=0, \, u_1^2+\ldots+u_n^2=1\}.$$ 
Denote by $R(\varphi,\Omega)$  the rotation in the plane $\{u_i=0, \, i=3,\ldots,n\}$ through an angle $\varphi$ about the origin $\Omega:$
$$u'_1=u_1\cos{\varphi}-u_2\sin{\varphi}, \quad u'_2=u_1\sin{\varphi}+u_2\cos{\varphi}, \quad u'_i=u_i, \; i=3, \ldots, n.$$

Let $$H_+=\{\vec u\in {\bf S }^{n-1}: u_2\ge 0\}, \quad H_-=\{\vec u\in {\bf S}^{n-1}: u_2\le 0\},$$ $$Q=\{\vec u\in {\bf S}^{n-1}: u_2=0,\; u_1>0\},\quad \bar Q=\{\vec u\in {\bf S}^{n-1}: u_2=0,\; 
u_1\ge 0\}.$$
Note that $H_-$ and $H_+$ are closed hemispheres of ${\bf S}^{n-1},\; \, \bar Q=Q\bigcup\Omega,\; $ and $\bar Q$ is
a hemisphere of the unit sphere $\Omega_2=\{\vec u\in {\bf S}^{n-1}: u_2=0\}$ bounded by $\Omega.$

\begin{lemma} Consider two points $y$ and $e_0$ in ${\bf S}^{n-1}.$ Suppose $y\in Q$ and 
$e_0\notin \bar Q.\\$ If $e_0\in H_+$, then any rotation $R(\varphi,\Omega)$ of $y$ with sufficiently small positive $\varphi$ decreases the distance between $y$ and $e_0.\\$
If $e_0\in H_-$, then any rotation $R(\varphi,\Omega)$ of $y$ with sufficiently small negative $\varphi$ decreases the distance between $y$ and $e_0.$
\end{lemma}
\begin{proof} Let $y$ be rotated into the point $y(\varphi)$. If the coordinate expressions of $y$ and $e_0$ are
$$y=(u_1,0,u_3,\ldots,u_n), \quad u_1>0; \qquad e_0=(v_1,v_2,\ldots,v_n), \; \, \mbox{then}$$ 
$$r(\varphi):=y(\varphi)\cdot e_0=u_1v_1\cos{\varphi}+u_1v_2\sin{\varphi}+u_3v_3+\ldots+u_nv_n.$$ 
Therefore, $\; r'(\varphi)=-u_1v_1\sin{\varphi}+u_1v_2\cos{\varphi},\; $ i.e. $\; r'(0)=u_1v_2. \;$ Then 
$$r'(0)>0 \quad \mbox{iff} \quad v_2>0, \quad \mbox{i.e.}\quad e_0\in \stackrel{\sf o}{H}_+;$$  
$$r'(0)<0 \quad \mbox{iff} \quad v_2<0  \quad \mbox{i.e.}\quad e_0\in \stackrel{\sf o}{H}_-.$$
That proves the lemma for $v_2\neq 0$. In the case $v_2=0$, by assumption ($e_0\notin \bar Q$) we have $v_1<0.$ In this case $r'(0)=0$, and $r''(0)=-u_1v_1>0$, i.e. $\varphi=0$ is a minimum point.
This completes the proof.
\end{proof}

\begin{prop} Let $Y$ be irreducible and $m=|Y| \ge n$. Suppose there are no closed great hemispheres $\bar Q$   in ${\bf S}^{n-1}$ such that $\bar Q$  contains $n-1$ points from $Y$ and $e_0$.  Then any vertex of $\Gamma_\psi(Y)$ has degree at least $n-1$.
\end{prop}
\begin{proof} Without loss of generality, we may assume that
$$\phi_{1,i}=\psi,\; \, i=2, \ldots, \deg{y_1}+1; \quad \phi_{1,i}>\psi, \; \, i=\deg{y_1}+2, \ldots,  m.$$

 Suppose $\deg{y_1}< n-1$. Then $\phi_{1,i}>\psi$ for $i=n, \ldots, m.$ Let us consider the great $(n-3)$ - dimensional 
sphere $\Omega$ in ${\bf S}^{n-1}$ that contains the points $y_2, \ldots, y_{n-1}.$ Then Lemma 3 implies that a rotation $R(\varphi,\Omega)$ of $y_1$ with sufficiently small $\varphi$ decreases $\theta_1$. It contradicts the irreducibility of $Y$.
\end{proof}

\begin{prop} If $Y$ is irreducible, $|Y|=n, \dim{\Delta_n}=n-1$, then $\deg{y_i} = n-1\; $ for all $\; i=1, \ldots, n.$ In other words, ${\Delta_n}$
is a regular simplex of edge lengths $\psi.$
\end{prop}
\begin{proof} Clearly, $\Delta_n$ is a spherical simplex. Denote by $F_i$ its facets, 
$$F_i:=\conv{\{y_1,\ldots,y_{i-1},y_{i+1},\ldots,y_n \}}.$$ 
Let for $\sigma\subset I_n:=\{1,\ldots,n\}$ 
$$F_\sigma:=\bigcap\limits_{i\in\sigma} {F_i}\,.$$ 

We claim for $ \; i\neq j \; $ that:
$$\mbox{\it If}\; \; e_0\notin F_{\{i,j\}}, \; \mbox{\it then \;} \phi_{i,j}=\psi. \eqno (5.4)$$ 

Conversely, from Lemma 3 follows that there exists a rotation $R(\varphi,\Omega_{ij})$ of $y_i$  (or $y_j$ if $e_0 \in F_i$) decreases $\theta_i$ (respectively, $\theta_j$), where $\Omega_{ij}$ is the great $(n-3)\, - \,$dimensional sphere contains $F_{\{i,j\}}$. It contradicts the irreducibility assumption for $Y$.

This yields, if there is no pair $\{i,j\}$ such that $e_0\in F_{\{i,j\}}$, then $\phi_{i,j}=\psi$ for all $i, j$ from $I_n$.

Suppose $e_0\in F_\sigma$, where $\sigma$ has maximal size and $|\sigma|>1$. Let $\bar\sigma=I_n\setminus\sigma$. From $(5.4)$ it follows that
$\phi_{i,j}=\psi\;$ if $\;i\in\bar\sigma\;$ or $\;j\in\bar\sigma.\; $ It remains to prove that $\phi_{i,j}=\psi\;$ for $\;i, j\in\sigma.\;$ 

Let $\Lambda$ be the intersection of the spheres of centers $y_i, \; i\in \bar\sigma,\; $ and radius $\psi$. Then $\Lambda$ is a sphere in ${\bf S}^{n-1}$ of dimension $|\sigma|-1.$ Note that $F_\sigma=$convex hull of $\{y_i: i\in\bar\sigma\}$, and for any fixed point $x$ from 
$F_\sigma$ (in particular for $x=e_0$) the distance $\dist(x,y)$ posses the same value (depending only on $x$) on the entire set $y\in\Lambda.$ Then  
$y_i, \; i\in\sigma,\; $ lie in $\Lambda$ at the same distance from $e_0$. It is clear that $Y$ is irreducible if and only if $y_i, \; i\in\sigma,\; $ in $\Lambda$ are vertices of a regular simplex of edge length $\psi.$

Finally, we have that all edges of $\Delta_n$ are of lengths $\psi$ as required.
\end{proof}

\begin{cor} If $n>3$, then $\Delta_4$ is a regular tetrahedron of edge lengths $\psi.$
\end{cor}
\begin{proof} 
Let us show that $\dim{\Delta_4}=3$.  In the converse case,  $\dim{\Delta_4}=2$, and from Theorem 4 follows that  $\Delta_4$ is a rhomb. Suppose $y_1y_3$ is the minimal length diagonal of 
$\Delta_4$.  Then $\phi_{2,4}>\psi$ (see Lemma 2).
Let us consider a sufficiently small turn of the facet $y_1y_2y_3$ round $y_1y_3$. If $e_0\notin y_1y_3$, then this turn decreases either $\theta_4$ (if $e_0\in y_1y_2y_3$) or $\theta_2$, a contradiction. In the case $e_0\in y_1y_3$ any turn of $y_2$ round $y_1y_3$ decreases $\phi_{2,4}$ and doesn't change $\theta_2$. Obviously, there is a turn such that $\phi_{2,4}$ becomes is equal to $\psi.$ That contradicts the irreducibility of $Y$ also.

\end{proof}

\medskip

\noindent{\bf 5-E. Irreducible sets in ${\bf S}^3$.}
\begin{lemma} If $Y\subset {\bf S}^3$ is irreducible, $|Y|=5$, then $\; \deg{y_i}\ge 3\; $
 for all $i$.
\end{lemma}
\begin{proof} ({\bf 1})
Let us show that $\dim{\Delta_5}=3$.  In the converse case,  $\dim{\Delta_5}=2$, and from Theorem 4 follows that  $\Delta_5$ is a convex equilateral pentagon. Suppose $y_1y_3$ is the minimal length diagonal of $\Delta_5$.  We have $\phi_{2,k}>\psi$ for $k>3.$
 Suppose $e_0\notin y_1y_3$. If $e_0\in y_1y_2y_3$ then any sufficiently small turn of the facet $y_1y_3y_4y_5$ round $y_1y_3$ decreases $\theta_4$ and $\theta_5$, otherwise it decreases $\theta_2$, a contradiction. 
In the case $e_0\in y_1y_3$ any turn of $y_2$ round $y_1y_3$ decreases $\phi_{2,k}$ for $k=4, 5$, and doesn't change $\theta_i$. It can be shown in the elementary way that there is a turn such that $\phi_{2,4}$ or $\phi_{2,5}$ becomes is equal to $\psi$, a contradiction.

In three dimensions there  exist only two combinatorial types of convex polytopes with 5 vertices: (A) and (B) (see Fig. 5). In the case (A) the arc $y_3y_5$ lies inside $\Delta_5$, and for (B): $y_2y_3y_4y_5$ is a facet of $\Delta_5.$

\medskip

\begin{center}
\begin{picture}(320,150)(-80,-75)
\put(57,-75){Fig. 5}

\put(-10,-65){(A)}

\put(-60,-20){\circle*{4}}

\put(20,-40){\circle*{4}}
\put(0,20){\circle*{4}}

\put(60,40){\circle*{4}}

\put(0,60){\circle*{4}}

\put(-60,-20){\line(4,-1){80}}
\put(-60,-20){\line(3,2){60}}
\put(-60,-20){\line(3,4){60}}

\put(0,20){\line(3,1){60}}
\put(20,-40){\line(1,2){40}}

\put(0,60){\line(0,-1){40}}
\put(0,60){\line(1,-5){20}}
\put(0,60){\line(3,-1){60}}
\put(20,-40){\line(-1,3){20}}
\thinlines
\multiput(-60,-20)(4,2){30}%
{\circle*{1}}

\put(-75,-21){$y_5$}
\put(27,-41){$y_2$}
\put(-13,24){$y_4$}
\put(64,44){$y_3$}
\put(3,65){$y_1$}


\put(170,-65){(B)}

\put(130,-20){\circle*{4}}

\put(190,-40){\circle*{4}}
\put(170,20){\circle*{4}}

\put(230,0){\circle*{4}}

\put(170,60){\circle*{4}}

\put(130,-20){\line(3,-1){60}}
\put(130,-20){\line(1,1){40}}
\put(130,-20){\line(1,2){40}}

\put(170,20){\line(3,-1){60}}
\put(190,-40){\line(1,1){40}}

\put(170,60){\line(0,-1){40}}
\put(170,60){\line(1,-5){20}}
\put(170,60){\line(1,-1){60}}
%

\put(115,-21){$y_5$}
\put(197,-41){$y_2$}
\put(157,24){$y_4$}
\put(234,4){$y_3$}
\put(173,65){$y_1$}

\end{picture}
\end{center}


\noindent ({\bf 2}) By $s_{ij}$ we denote the arc $y_iy_j$, and by $s_{ijk}$ denote the triangle $y_iy_jy_k.$  Let $\tilde s_{ijk}$ be the intersection of the great $2-$hemisphere $Q_{ijk}$ and $\Delta_5$, where $Q_{ijk}$ contains $y_i, y_j, y_k$ and bounded by the great circle passes through $y_i, y_j$.
Proposition 2 yields: if there are no $i, j, k$ such that $e_0\in \tilde s_{ijk}$, then 
$\deg{y_i}\ge 3$ for all $i$. 

It remains to consider all cases $e_0\in \tilde s_{ijk}$. Note that
for (A) $\tilde s_{ijk}\ne s_{ijk}$ only for three cases: $i=1, 2,4;$ where $j=3,\; k=5,$ or $j=5,\; k=3$ ($\tilde s_{i35}=\tilde s_{i53}$). 

\medskip

\noindent ({\bf 3})  Lemma 1 yields that $\deg{y_k}>0$. Now we consider the cases $\deg{y_k}=1, 2$.

\medskip

 \centerline{\em If $\; \deg{y_k}=1,\; \phi_{k,\ell}=\psi,\; $ then $\; e_0\in s_{k\ell}.$} 

\medskip

Indeed, otherwise
there exists the great circle
$\Omega$ in ${\bf S}^3$ such  that $\Omega$ contains $y_\ell$, and  the great sphere passes through $\Omega$ and $y_k$ doesn't pass through $e_0$. Then Lemma 3 implies that a rotation $R(\varphi,\Omega)$ of $y_k$ with sufficiently small $\varphi$ decreases $\theta_k$ - a contradiction.
 
Since $\theta_0<\psi$, $e_0$ can not be a vertex of $\Delta_5.$ Therefore, $e_0$ lies inside $s_{k\ell}$. From this follows if $s_{ij}$ for any $j$ doesn't intersect $s_{k\ell}$, then $\deg{y_i}\ge 2.$

Arguing as above it is easy to prove that 

\medskip

\centerline{\em If $\; \deg{y_k}=2,\; \phi_{k,i}=\phi_{k,j}=\psi,\; $ then $\; e_0\in \tilde s_{ijk}.$}

\medskip

\noindent ({\bf 4}) Now we prove that $\deg{y_k}\ge 2$ for all $k.$ 
Conversely, $\deg{y_k}=1, \; e_0\in s_{k\ell}.$ 

a). First we consider the case when $s_{k\ell}$ is an ``external" edge of $\Delta_5$. For the type (A) that means $s_{k\ell}$ differs from $s_{35}$, and for (B) it is not $s_{35}$ or $s_{24}$. Since $\Delta_5$ is convex, there exists the great $2-$sphere $\Omega_2$ passes through $y_k, y_\ell$ such that 3 other points $y_i, y_j, y_q$ lie inside the hemisphere $H_+$ bounded by $\Omega_2.$ Let $\Omega$ be the great circle in $\Omega_2$ that contains $y_\ell$ and is orthogonal to the arc $s_{k\ell}$. Then (Lemma 3) there exists a small turn of $y_i, y_j, y_q$ round $\Omega$ that simultaneously  decreases $\theta_i, \theta_j, \theta_q$ - a contradiction.

b). For the type (A) when $\deg{y_3}=1, \; \phi_{3,5}=\psi, \; e_0\in s_{35}$; we claim that $s_{124}$ is a regular triangle with side length $\psi.$ Indeed, from a) follows that $\deg{y_i}\ge 2$ for $i=1, 2, 4.$ Moreover, if $\deg{y_i}=2$, then $e_0=s_{35}\bigcap s_{124}.$ Therefore, in any case, $\phi_{1,2}=\phi_{1,4}=\phi_{2,4}=\psi.$ 
We have the arc  $s_{35}$ and the regular triangle $s_{124}$, both are with edge lengths $\psi$. Then from Lemma 2 follows that some $\phi_{i,j}<\psi$ - a contradiction.

c). Now for the type (B) consider the case:  $\deg{y_3}=1, \; \phi_{3,5}=\psi, \; e_0\in s_{35}$. Then 
 for $y_2$ we have: $\deg{y_2}=1$ only if $\phi_{2,4}=\psi, \; e_0=s_{24}\bigcap s_{35};\; $ $\deg{y_2}=2$ only if $\phi_{2,4}=\phi_{2,5}=\psi$; and $\phi_{2,4}=\phi_{1,2}=\phi_{2,5}=\psi$ if $\deg{y_2}=3$.
Thus, in any case, $\phi_{2,4}=\psi.$ We have two intersecting diagonals $s_{24}, s_{35}$ of lengths $\psi.$ Then Lemma 2 contradicts the assumption that $Y$ is a $\psi$-code. This contradiction concludes the proof that $\deg{y_k}\ge 2$ for all $k$.

\medskip

\noindent ({\bf 5}) Finally let us prove that $\deg{y_k}\ge 3$ for all $k.$ Assume the converse. Then $\; \deg{y_k}=2, \; e_0\in \tilde s_{ijk},\; $ where $\; \phi_{k,i}=\phi_{k,j}=\psi.$ 

{\bf Case facet}{\bf:} Let $s_{ijk}$ be a facet of $\Delta_5,$ and $e_0\notin s_{ij}$. By the same argument as in ({\bf 4}a), where $\Omega_2$ be the great sphere  contains $s_{ijk}$, and $\Omega$ be the great circle passes through $y_i, y_j,$ we can prove that there exists a shift decreases $\theta_\ell, \theta_q$ for two other points $y_\ell, y_q$ from $Y$, a contradiction.

If $e_0\in s_{ij}$, then any turn of $s_{{\ell}q}$ round $\Omega$ doesn't change $\theta_\ell$ and $\theta_q$. However, if this turn is in a positive direction, then it decreases $\phi_{k,\ell}$ and $\phi_{k,q}$. Clearly, there exists a turn when $\phi_{k,\ell}$ or $\phi_{k,q}$ is equal to $\psi$ - a contradiction.

It remains to consider all cases where $s_{ijk}$ is not a facet. Namely, there are the following cases: 
$s_{124},\; s_{135}$ (type (A) and type (B)), $\; s_{234}$ (type (B)).

{\bf Case $s_{124}$:} We have $\deg{y_1}=2,\; \phi_{1,2}=\phi_{1,4}=\psi,\; e_0\in s_{124}.$ Consider a small turn of $y_3$ round $s_{24}$ towards $y_1$. If $e_0\notin s_{24}$, then this turn decreases $\theta_3.$ Therefore, the irreducibility yields $\phi_{3,5}=\psi.$ In the case $e_0\in s_{24},\; \theta_3'=\theta_3,$ but $\phi_{1,3}$ decreases. It again implies  $\phi_{3,5}=\psi.$ Since $s_{35}$ cannot intersects a regular triangle $s_{124}$ [see Lemma 2, ({\bf 4}b)], $\phi_{2,4}>\psi.$
Then $\deg{y_2}=\deg{y_4}=3.$ (Since $ e_0\in s_{124},\; \deg{y_2}=2$ only if   $\phi_{2,4}=\psi.$)
Thus we have three isosceles triangles $s_{243}, s_{241}, s_{245}$. Using this and $\phi_{3,5}=\psi,$ we obviously have $\phi_{1,i}<\psi$ for $i=3, 5,$ - a contradiction.

{\bf Case $s_{135}$}(type (B)) is equivalent to the {\bf Case $s_{124}$}.

{\bf Case $s_{135}$}(type (A)){\bf:} This case has two subcases: $\tilde s_{351},\; \tilde s_{153}$.

In the subcase $\tilde s_{135}$ we have 
$\deg{y_1}=2, \; \phi_{1,3}=\phi_{1,5}=\psi, \; e_0\in \tilde s_{135}.\\$ If $e_0\notin s_{135}$, then any turn of $y_1$ round $s_{35}$ decreases $\theta_1$ (Lemma 3). Then $e_0\in s_{135}$. 
Clearly, any small turn of $y_2$ round $s_{35}$ increases $\phi_{2,4}.$ On the other hand, this turn decreases $\theta_2$ (if $e_0\notin s_{35}$) and $\phi_{1,2}$. Arguing as above, we get a contradiction. 
The subcase $\tilde s_{315}$, where $\phi_{3,5}=\psi$, can be proven by the same arguments as Case $s_{124}$. 

{\bf Case $s_{234}$}(type (B)){\bf:} This case has two subcases: $\tilde s_{243},\; \tilde s_{234}$.

It is not hard to see that  $\tilde s_{243}$ follows from the {\bf Case facet}, and $\tilde s_{234}$ can be proven in the same way as the subcase $\tilde s_{135}$. This concludes the proof.
\end{proof}

Lemma 4 yields that the degree of any vertex of $\Gamma_\psi(Y)$ is not less than 3. This implies that at least one vertex of $\Gamma_\psi(Y)$ has degree 4. Indeed, if all vertices of $\Gamma_\psi(Y)$ are of degree 3, then the sum of the degrees equals 15, i.e. is not an even number. 
There exists only one type of $\Gamma_\psi(Y)$ with these conditions (Fig. 6). The lengths of all edges of $\Delta_5$ except $y_2y_4$, $y_3y_5$ are equal to $\psi$.
For fixed $\phi_{2,4}=\alpha, \; \Delta_5$ is uniquely defined up to isometry. Therefore, we have the 1-parametric family $P_5(\alpha)$ on ${\bf S}^3.\;$ If $\phi_{3,5}\ge\phi_{2,4}$, then $z\ge\cos{\alpha}\ge 2z-1.$

\begin{center}
\begin{picture}(320,140)(-80,-70)
\put(57,-65){Fig. 6: $P_5(\alpha)$}

\put(10,-20){\circle*{5}}

\put(90,-40){\circle*{5}}
\put(70,20){\circle*{5}}

\put(130,40){\circle*{5}}

\put(70,60){\circle*{5}}

\thicklines
\put(10,-20){\line(4,-1){80}}
\put(10,-20){\line(3,2){60}}

\put(70,20){\line(3,1){60}}
\put(90,-40){\line(1,2){40}}

\put(70,60){\line(-3,-4){60}}
\put(70,60){\line(0,-1){40}}
\put(70,60){\line(1,-5){20}}
\put(70,60){\line(3,-1){60}}

\thinlines
\multiput(90,-40)(-1,3){20}%
{\circle*{1}}

\put(70,-14){$\alpha$}

\put(-5,-21){$y_5$}
\put(97,-41){$y_2$}
\put(57,24){$y_4$}
\put(134,44){$y_3$}
\put(73,65){$y_1$}
 
\end{picture}
\end{center}


Thus Theorem 4, Corollary 4 and Lemma 4 for $n=4$ yield:

\begin{theorem} Let $Y\subset {\bf S}^3$ be an irreducible set, $\; |Y|=m\le 5.\; $ Then $\; \Delta_m \; $ 
 for $\; 2\le m\le 4$ is  a regular simplex of edge lengths $\psi$, and $\; \Delta_5$ is isometric to $P_5(\alpha)$ for some
$\alpha\in [\psi,\arccos{(2z-1)}].$
\end{theorem}

\medskip

\noindent{\bf 5-F. Optimization problem.} We see that if $Y$ is optimal, then for some cases 
it can be defined up to isometry. 
For fixed $ y_i\in {\bf S}^{n-1},\; i=1,\ldots,m;\; $ the function $H$ depends only on a position $y=-y_0=e_0\in {\bf S}^{n-1}.$ Let
$$H_m(y):=f(1)+f(-y\cdot y_1)+\ldots+f(-y\cdot y_m),$$ 
i.e. $H_m(y)=H(-y;Y).\; $ 

Thus for $h_m$ we have the following $(n-1)$-dimensional optimization problem: 

$$h_m=\max\limits_y{\{H_m(y)\}}$$ {\em subjects to the constraint}
$$y\in T(Y,\theta_0):=\{y\in\Delta_m\subset {\bf S}^{n-1}:\; y\cdot y_i\ge t_0,\; i=1,\ldots, m\}.  $$ 

We present an efficient numerical method for this problem in the next section.

\section {On calculations of $h_m$}
In this technical section we explain how to find an upper bound on  $\; h_m\; $ for $\; n=4, \; m\le 6$. 
Note that Theorem 5 gets for computation of $h_m$ a low-dimensional optimization problem 
(see {\bf 5-F}). Our first approach for this problem was to apply numerical methods \cite{Mus}. However, that is a nonconvex constrained optimization problem. In this case, the Nelder-Mead simplex method and other local improvements methods cannot guarantee finding a global optimum. It's possible (using estimations of derivatives) to organize computational process in such way that it gives a global optimum. However, such solutions are very hard to verify and some mathematicians don't accept that kind of proofs. Fortunately, using geometric approach, estimations of $h_m$ can be reduced to relatively simple computations.

Throughout this section we use the function $\tilde f(\theta)$ defined for 
$f\in {\it\Phi}^*(z)$ by 
$$ \tilde f(\theta):=\left\{
\begin{array}{l}
f(-\cos{\theta}) \quad  0\le\theta\le\theta_0=\arccos{t_0} \; \mbox{ (see Definition 4)}\\
-\infty \qquad \quad \; \;  \theta>\theta_0
\end{array} 
\right.
$$
Since $f\in {\it\Phi}^*(z)$,  $\tilde f(\theta)$ is a monotone decreasing function in $\theta$ on $[0,\theta_0]$.

\medskip

\noindent{\bf 6-A. The case m=2.} 
Suppose $m=2$  and $Y$ is optimal for $f\in {\it\Phi}^*(z)$. Then $\Delta_2=y_1y_2$ is an arc of length $\psi, \; e_0\in\Delta_2,$ and $\; \theta_1+\theta_2=\psi,$ where $\theta_i\le\theta_0$ (see Lemma 1 and $(5.2)$). The efficient function $F(\theta_1,\theta_2)=f(1)+
\tilde f(\theta_1)+\tilde f(\theta_2)$  
is a symmetric function in $\theta_1, \theta_2.$ 

 We can assume
that $\theta_1\le \theta_2$, then $\theta_1\in [\psi-\theta_0,\psi/2].$
Since $\Theta_2(\theta_1):=\psi-\theta_1$ is a monotone decreasing function, $\tilde f(\Theta_2(\theta_1))$ is a monotone increasing function in $\theta_1.$
Thus for any $\theta_1\in [u,v]\subset [\psi-\theta_0, \psi/2]$ we have
$$ F(\theta_1,\theta_2)\le \Phi_2([u,v]):=f(1)+\tilde f(u)+\tilde f(\psi-v). $$

Let $\; u_1=\psi-\theta_0,\; u_2,\, \ldots,\, u_{N},\; u_{N+1}=\psi/2\; $ be points in $\; [\psi-\theta_0,\psi/2]\; $ such that $\; u_{i+1}=u_i+\varepsilon,\; $ where $\; \varepsilon=(\theta_0-\psi/2)/N.\; $ 
If $\; \theta_1\in [u_i,u_{i+1}],\; $ then $\; h_2=H(y_0;Y)=F(\theta_1,\theta_2)\le 
\Phi_2([u_i,u_{i+1}]). \; $ Thus
$$h_2\le \lambda_2(N,\psi,\theta_0):=\max\limits_{1\le i\le N} 
\{\Phi_2(s_i)\}, \; \mbox{ where } \; s_i:=[u_i,u_{i+1}].\; $$
Clearly, $\lambda_2(N,\psi,\theta_0)$ tends to $\, h_2\, $ as $\; N\to\infty\; $ ($\varepsilon\to 0$). 

That implies a very simple method for calculation of $h_2.$ Now we extend this approach to higher $m$.

\medskip

\noindent{\bf 6-B. The function $\Theta_k$.}
Suppose we know (up to isometry) optimal $Y=\{y_1, \ldots, y_m\}\subset {\bf S}^{n-1}$. Let us assume that 
 $\dim{\Delta_m}=n-1$, and $\; V:= $ convex hull of $\{y_1\ldots y_{n-1}\} $ is a facet of $\Delta_m.$ 
Then  $\; \mbox{rank}{\{y_1, \ldots, y_{n-1}\}}=n-1,\; $ and $Y$ belongs to the hemisphere $H_+$, where $H_+$ contains $Y$ and bounded by the great sphere  $\tilde S$  passes through $V$.

Let us show that any $y=y_+\in H_+$ is uniquely determined by the set of distances $\theta_i=\dist(y,y_i),\; i=1, \ldots, n-1.$ Indeed, there are at most  two solutions: $\; y_+\in H_+\; $ and $\; y_-\in H_-\; $ of the quadratic equation 
$$y\cdot y=1 \; \mbox{ with } \; y\cdot y_i=\cos{\theta_i}, \; i=1,\ldots,n-1. \eqno (6.1)$$
Note that $y_+=y_-$ if and only if $y\in \tilde S.$  

This implies that $\theta_k,\; k\ge n $ is determined by $\; \theta_i,\; i=1,\ldots, n-1;$
$$ \theta_k=\Theta_k(\theta_1,\ldots,\theta_{n-1}).$$
It is not hard to solve $(6.1)$ and, therefore, to give an explicit expression for $\Theta_k.$

For instance, let $\Delta_n$ be a regular simplex of edge lengths $\pi/3$. (We need this case for $n=3,4$.) Then 
\footnote{I am very grateful to referees  for these explicit formulas.}
$$
\cos\theta_3=\cos\Theta_3(\theta_1,\theta_2)=
$$
$$
=\frac{1}{3}\left(\cos\theta_1+\cos\theta_2+\sqrt{6-8[\cos\theta_1\cos\theta_2+
(\cos\theta_2-\cos\theta_1)^2]}\right);
$$
$$
\cos\theta_4=\cos\Theta_4(\theta_1,\theta_2,\theta_3)=\frac{1}{4}\Bigl(\cos\theta_1+\cos\theta_2+\cos\theta_3+
$$
$$
\sqrt{10}\sqrt{1+\cos\theta_1\cos\theta_2+\cos\theta_1\cos\theta_3+\cos\theta_2\cos\theta_3
-\frac{3}{2}(\cos^2\theta_1+\cos^2\theta_2+\cos^2\theta_3)}\Bigr).
$$

\noindent{\bf 6-C. Extremal points of $\Theta_k$ on $D$.}
Let $\; {\bf{a}}=(a_1,\ldots,a_{n-1}),\; $ where $\; 0<a_i\le \theta_0<\psi.$  (Recall that $ \phi_{i,j}=\dist(y_i,y_j);\; \cos{\psi}=z; \; \cos{\theta_0}=t_0.$) Now we consider a domain
$D({\bf{a}})$ in $H_+$, where 
$$D({\bf{a}})=\{y\in H_+: \; \dist(y,y_i)\le a_i,\; \, 1\le i\le n-1\}.$$
In other words,
$D({\bf{a}})$ is the intersection of the closed caps
$\cp(y_i,a_i)$ in $H_+$:
$$D({\bf{a}})=\bigcap\limits_{i=1}\limits^{n-1} {\cp(y_i,a_i)}\bigcap H_+ . $$

Suppose $\dim{D({\bf{a}})}=n-1$. Then $D({\bf{a}})$ has ``vertices", ``edges", and ``$k$-faces" for $k\le n-1.$
Indeed, let $$\sigma\subset I:=\{1, \ldots, n-1\},\quad 0<|\sigma|\le n-1;$$ 
$$\tilde F_\sigma:=\{y\in D({\bf{a}}):\; \dist(y,y_i)=a_i \; \, \forall \; i\in\sigma\}.$$
It is easy to prove that $\dim{\tilde F_\sigma}=n-1-|\sigma|;\; \tilde F_\sigma$ belongs to the boundary $B$ of $D({\bf{a}})$; and if $\sigma\subset\sigma'$, then  $\tilde F_{\sigma'}\subset \tilde  F_\sigma$. 
Actually, $D({\bf{a}})$ is combinatorially equivalent to an $(n-1)$-dimensional simplex.

Now we consider the minimum of $\Theta_k(\theta_1,\ldots,\theta_{n-1})$ on $D({\bf{a}})$ for $k\ge n$. In other words, we are looking for a point $p_k({\bf{a}})\in D({\bf{a}})$ such that
$$\dist(y_k,p_k({\bf{a}}))=\dist(y_k,D({\bf{a}})).$$ 
Since $\; \phi_{i,k}\ge \psi>\theta_0,\; $ all $y_k$ lie outside  $D({\bf{a}})$. Clearly, $\Theta_k$ achieves its minimum at some point in $B$. Therefore, there is $\sigma\subset I$ such that $$p_k({\bf{a}})\in \tilde F_\sigma \eqno (6.2)$$ 

Suppose $\sigma=I,$ then $\tilde F_\sigma$ is a vertex of $D({\bf{a}})$. Let us denote this point by $p_*({\bf{a}})$. Note that the function $\Theta_k$ at the point $p_*({\bf{a}})$ is equal to $\Theta_k({\bf{a}})$.

Let $\sigma_k({\bf{a}})$ denote $\sigma\subset I$ of the maximal size such that $\sigma$ satisfies $(6.2)$.
Then for $\sigma_k({\bf{a}})=I,\; p_k({\bf{a}})=p_*({\bf{a}})$, and for $|\sigma_k({\bf{a}})|<n-1, \; p_k({\bf{a}})$ belongs to the open part of $\tilde F_{\sigma_k({\bf{a}})}$.

Consider $n=3$. There are two cases for $p_k({\bf{a}})$ (see Fig. 7): $\; p_3({\bf{a}})=p_*({\bf{a}})=
\tilde F_{\{1,2\}},$ and $p_4({\bf{a}})$ is the intersection in $H_+$ of the great circle passes through $y_1, \; y_4$, and the circle $\tilde S(y_1,a_1)$ of center $y_1$ and radius $a_1$ ($\tilde F_{\{1\}}\subset \tilde S(y_1,a_1)$).
The same holds for all dimensions.

Denote by $S_\sigma(k)$ the great $|\sigma|-$dimensional sphere passes through $y_i,\; i\in\sigma,$ and $y_k.$ Let $\tilde S(y_i,a_i)$ be the sphere of center $y_i$ and radius $a_i$; and for $\sigma\subset I$
$$\tilde S_\sigma:=\bigcap\limits_{i\in\sigma} \tilde S(y_i,a_i).$$ 
Denote by $s(\sigma,k)$ the intersection of $S_\sigma(k)$ and $\tilde S_\sigma$ in $H_+$,
$$s(\sigma,k)=S_\sigma(k)\bigcap  \tilde S_\sigma\bigcap H_+$$

\begin{center}
\begin{picture}(320,140)(-160,-70)

\put(-150,-30){\circle*{4}}
\put(-60,-30){\circle*{4}}
\put(-130,45){\circle*{4}}
\put(-50,45){\circle*{4}}
\put(-120,10){\circle*{3}}
\put(-109.5,0){\circle*{3}}

\qbezier (-100,-30) (-100,-22) (-102,-15)
\qbezier (-102,-15) (-105,-8) (-109,-1)
\qbezier (-109,-1) (-113,5) (-120,10)

\qbezier (-120,10) (-125,4) (-128,-2)
\qbezier (-128,-2) (-131,-8) (-132,-22)

\put(-132,-22){\line(0,-1){8}}

\put(-150,-30){\line(1,0){90}}

\multiput(-150,-30)(5,3.75){20}%
{\circle*{1}}

\put(-155,-40){$y_1$}
\put(-65,-40){$y_2$}
\put(-144,40){$y_3$}
\put(-45,40){$y_4$}
\put(-127,-23){$D({\bf{a}})$}
\put(-170,10){$H_+$}
\put (-99,-20){$\tilde F_{\{1\}}$}
\put (-150,-4){$\tilde F_{\{2\}}$}
\put (-130,15){$p_*({\bf{a}})$}
\put (-104,-3){$p_4({\bf{a}})$}

\put(-118,-65){Fig. 7}

\put(80,-65){Fig. 8}

\qbezier (60,-30) (100,-10) (140,-10)
\qbezier (40,20) (80,40) (120,40)

\qbezier (77,-37) (72,3) (47,43)
\qbezier (130,-25) (125,15) (100,55)

\put(76,5){$E({\bf b},{\bf{a}})$}
\put(27,-36){$\theta_2={ b}_2$}
\put(10,11){$\theta_2=a_2$}
\put(63,-45){$\theta_1={ b}_1$}
\put(115,-33){$\theta_1=a_1$}

\end{picture}
\end{center}

\begin{lemma} Suppose $D({\bf{a}})\ne\emptyset, \; 0<a_i\le \theta_0\; $ for all $i$, and
$\; k\ge n.\; $ 
Then $(i) \; \; p_k({\bf{a}})\in s(\sigma_k({\bf{a}}),k),\\(ii) $
if  $\; s(\sigma,k)\ne\emptyset, \; |\sigma|<n-1$, then $s(\sigma,k)$ consists of the one point $p_k({\bf{a}}).$
\end{lemma}
\begin{proof} 
$(i)$ Let $\theta_k^*:=\Theta_k(p_k({\bf{a}}))=\dist(y_k,p_k({\bf{a}})).$ Since $\Theta_k$ achieves its minimum at $p_k({\bf{a}})$, the sphere $\tilde S(y_k,\theta_k^*)$ touches the sphere 
$\tilde S_{\sigma({\bf{a}})}$ at $p_k({\bf{a}})$. If some sphere touches the intersections of spheres,
then the touching point belongs to the great sphere passes through the centers of these spheres. Thus 
$\; p_k({\bf{a}})\in S_{\sigma({\bf{a}})}(k)$.

\noindent $(ii)\; $ Note that $\, s(\sigma,k)\, $ belongs to the intersection in $\, H_+\, $ of the spheres  
$\; S(y_i,a_i), \\ i\in\sigma,$ and $S_\sigma(k)$. 
Any intersection of spheres is also a sphere. Since
$$\dim{S_\sigma(k)}+\dim{\tilde S_\sigma}=n-1,$$
this intersection is empty, or is a $0-$dimensional sphere (i.e. 2-points set). In the last case, one point lies in $H_+$, and another one in $H_-.\; $ Therefore, $s(\sigma,k)=\emptyset,$ or 
$s(\sigma,k)=\{p\}.$ Denote by $\sigma'$ the maximal size $\sigma'\supset\sigma$ such that
$s(\sigma',k)=\{ p\}.$ It is not hard to see that $\tilde S(y_k,\dist(y_k,p))$ touches
$\tilde S_{\sigma'}$ at $p$. Thus $\;  p=p_k({\bf{a}})$.
\end{proof}

Lemma 5 implies a simple method for calculations of the minimum of $\Theta_k$ on $D({\bf{a}})$. For this we can consider $s(\sigma,k), \; \sigma\subset I$, and if $s(\sigma,k)\ne\emptyset,$ then
$s(\sigma,k)=\{p_k({\bf{a}})\},$ so then $\Theta_k$ attains its minimum at this point. In the case 
when $\Delta_n$ is a  simplex we can find the minimum by very simple method.

\begin{cor} Suppose $|Y|=n, \; 0<a_i\le \theta_0\; $ for all $i$, and $D({\bf{a}})$ lies inside $\Delta_n.$  Then 
$$\theta_n\ge\Theta_n(a_1,\ldots,a_{n-1}) \; \mbox{ for all }\; y\in D({\bf{a}}).$$
\end{cor}
\begin{proof} Clearly, $\Delta_n$ is a simplex. Since  $D({\bf{a}})$ lies inside $\Delta_n$, for $|\sigma|<n-1$ the intersection of $\tilde S_\sigma$ and $S_\sigma(k)$ is empty. Thus $p_n({\bf{a}})=p_*({\bf{a}}).$  
\end{proof}

\medskip

\noindent{\bf 6-D. Upper bounds on $H_m$.} 
Suppose $\; \dim{\Delta_m}=n-1, \; $ and $\; y_1\ldots y_{n-1}$ is a facet of $\Delta_m$.  
Then (see {\bf 5-F} for the definition of $H_m$ and $T(Y,\theta_0)$)
$$ H_m(y)=F(\theta_1,\ldots,\theta_{n-1},\Theta_n,\ldots,\Theta_m)=\tilde F_m(\theta_1,\ldots,\theta_{n-1}),$$
where
$$\tilde F_m(\theta_1,\ldots,\theta_{n-1}):=f(1)+\tilde f(\theta_1)+\ldots+\tilde f(\theta_{n-1})+\tilde f(\Theta_n(\theta_1,\ldots,\theta_{n-1}))
+\ldots$$ $+\, \tilde f(\Theta_m(\theta_1,\ldots,\theta_{n-1})).$

\begin{lemma} Suppose $f\in{\it\Phi}^*(z), \; 
  |Y|=m,\; \dim{\Delta_m}=n-1,\; y_1\ldots y_{n-1}$ is a facet of $\Delta_m, \; \dist(y_i,y_j)\ge \psi > \theta_0\; $ for $\; i\ne j$,
$\; 0\le{ b}_i<a_i\le\theta_0\; $ for $\; i=1,\ldots,n-1;$
and $\; \Theta_k({\bf{a}})\le\theta_0\; $ for all $\; k\ge n.\; $
If $\; D({\bf{a}})\ne\emptyset,\; $ then 
$$H_m(y)\le \Phi_Y({\bf b},{\bf{a}}) \quad \mbox{for any} \quad y\in E({\bf b},{\bf{a}}):=D({\bf{a}})\setminus U({\bf b}),$$
where 
$$\Phi_Y({\bf b},{\bf{a}}):=f(1)+\tilde f({ b}_1)+\ldots+\tilde f({ b}_{n-1})+
\tilde f(\Theta_n(p_n({\bf{a}})))+\ldots+\tilde f(\Theta_m(p_m({\bf{a}}))),$$
$$U({\bf b}):=\bigcup\limits_{i=1}^{n-1} \cp(y_i,{b}_i).$$
\end{lemma}
\begin{proof} We have for $1\le i\le n-1$ and $y\in E({\bf b},{\bf{a}}):\; \theta_i\ge { b}_i$ (Fig. 8). By the monotonicity assumption  this implies 
$\tilde f(\theta_i)\le \tilde f({ b}_i).$ On the other hand, 
$\; y\in D({\bf{a}}).\; $     Then Lemma 5 yields 
$\; \tilde f(\theta_k)\le \tilde f(\Theta_k(p_k({\bf{a}}))))\; $ for $\; k\ge n.$
\end{proof}

From Corollary 5 and Lemma 6 follow
\begin{cor} Let $|Y|=n$. Suppose $f, \; {\bf a},\; {\bf b},$ and $Y$ satisfy the assumptions of Lemma 6 and Corollary 5. Then for any $\; y\in E({\bf b},{\bf{a}}):$
$$H_m(y)\le f(1)+\tilde f({ b}_1)+\ldots+\tilde f({ b}_{n-1})+
\tilde f(\Theta_n({\bf{a}})).$$
\end{cor}

Let $K(n,\theta_0):=[0,\theta_0]^{n-1},\; $ i.e $\; K(n,\theta_0)$ is an $(n-1)-$dimensional  cube of side length $\theta_0$. Consider for $K(n,\theta_0)$ the cubic grid $L(N)$ of sidelength $\varepsilon,$ where $\varepsilon=\theta_0/N$ for given positive integer $N$. Then the grid (tessellation) $L(N)$ consists of $N^{n-1}$ cells, 
any cell $c\in L(N)$ 
is an $(n-1)-$dimensional cube of sidelength $\varepsilon,$ and for any point $(\theta_1,\ldots,\theta_{n-1})$ in $c$ we have
$${ b}_i(c)\le\theta_i\le a_i(c),  \quad a_i(c)={ b}_i(c)+\varepsilon, \quad
 i=1,\ldots,n-1.$$
Let $\tilde L(N)$ be the subset of cells $c$ in $L(N)$ such that 
$D({\bf{a}}(c))\ne\emptyset.$ There exists $c\in L(N)$ such that $H_m$ attains its maximum on $T(Y,\theta_0)$ at some point in $E({\bf b}(c),{\bf{a}}(c))$. Therefore, Lemma 6 yields

\begin{lemma} Suppose $f$ and $Y$ satisfy the assumptions of Lemma 6, $N$ is a positive integer, and $y\in\Delta_m \;$ is such that $\dist(y,y_i)\le \theta_0$ for all $i$. Then
$$H_m(y)\le \max\limits_{c\in\tilde L(N)}  \{\Phi_Y({\bf b}(c),{\bf{a}}(c))\}$$
\end{lemma}

\medskip

\noindent{\bf 6-E. Upper bounds on $h_m$.}
Suppose $\Delta_m$ is a regular simplex of edge length $\psi.$
Then the efficient function $F$ is a symmetric function in the variables $\theta_1, \ldots, \theta_m$. Consider this problem only on the domain 
$$\Lambda:=\{y\in\Delta_m:\; \psi-\theta_0\le\theta_1\le\theta_2\le \ldots\le\theta_m\le \theta_0\}.$$

Let $L_\Lambda(N)$ be the subset of cells $c$ in $\tilde L(N)$ such that 
$c\bigcap\Lambda\ne\emptyset.$ If $c\in L_\Lambda(N)$ is such that $E({\bf b}(c),{\bf{a}}(c)))$ lies inside $\Delta_m$,\footnote{Clearly, it holds for all $c\in L_\Lambda(N)$ if $N$ is sufficiently large.} then 
we have an explicit expression for $\Phi_m(c):=\Phi_Y({\bf b}(c),{\bf{a}}(c))$ (see Corollary 6).  For $n=4$, Theorem 5 implies that $\Delta_m$ is a regular simplex, where $m=2, 3, 4.$   Thus from Lemma 7 follows
$$h_m\le \lambda_m(N,\psi,\theta_0):=\max\limits_{c\in L_\Lambda(N)} 
\{\Phi_m(c)\}.$$

Now we consider the case $n=4,\; m=5.$ Theorem 5 yields: $\Delta_5$ is isometric to $P_5(\alpha)$ for some $\alpha\in [\psi,\psi':=\arccos{(2z-1)}]$ (see Fig. 6). Let the vertices $y_1, y_2, y_3$ of
$P_5(\alpha)$ be fixed. Then the vertices $y_4(\alpha),\; y_5(\alpha)$ are uniquely determined by $\alpha.$

Note that for any $y\in D(\theta_0,\theta_0,\theta_0)$ the distance 
$\theta_4(\alpha):=\dist(y,y_4(\alpha))$ increases, and $\theta_5(\alpha)$ decreases whenever $\alpha$ increases.
Let $\alpha_1=\psi,\, \alpha_2,\ldots, \alpha_N, \, \alpha_{N+1}=\psi'$ be points in $[\psi,\psi']$ such that
$\alpha_{i+1}=\alpha_i+\epsilon$, where $\epsilon=(\psi'-\psi)/N.$
 Then 
$$\theta_4(\alpha_i)<\theta_4(\alpha_{i+1}), \quad \theta_5(\alpha_i)>\theta_5(\alpha_{i+1}),$$
so then
$$\tilde f(\theta_4(\alpha_i))>\tilde f(\theta_4(\alpha_{i+1})), \quad \tilde f(\theta_5(\alpha_i))<
\tilde f(\theta_5(\alpha_{i+1})).$$
Combining this with Lemma 7, we get
$$h_5\le \lambda_5(N,\psi,\theta_0):=f(1)+ 
{\max\limits_{c\in \tilde L(N)}\{R_{1,2,3}(c)+\max\limits_{1\le i\le N}\{R_{4,5}(c,i)\}\}},$$
$$R_{1,2,3}(c)=\tilde f({ b}_1(c))+\tilde f({ b}_2(c))+\tilde f({ b}_3(c)),$$ $$
R_{4,5}(c,i)=\tilde f(\Theta_4(p_4({\bf{a}}(c),\alpha_{i})))+\tilde f(\Theta_5(p_5({\bf{a}}(c),\alpha_{i+1}))),  $$
where $p_k({\bf{a}},\alpha)=p_k({\bf{a}})$ with $y_k=y_k(\alpha).$
  
Clearly, $\lambda_m(N+1,\psi,\theta_0)\le\lambda_m(N,\psi,\theta_0)$. It's not hard to  show  that
$$h_m=\lambda_m(\psi,\theta_0):=\lim_{N\to\infty}{\lambda_m(N,\psi,\theta_0)}.$$

Finally let us consider the case: $n=4,\; m=6.$ In this case, we give an upper bound on $h_6$ by separate argument.
\begin{lemma}
 Let $\; n=4, \; f\in{\it\Phi}^*(z), \; \sqrt{z}>t_0>z,\; $ 
$ \theta_0'\in [\arccos{\sqrt{z}},\theta_0].$ Then
 $$h_6\, \le\, \max{\{\; \tilde f(\theta_0')+\lambda_5(\psi,\theta_0),\; 
f(-\sqrt{z})+\lambda_5(\psi,\theta_0')\; \}}.$$
\end{lemma}
\begin{proof}
Let $Y=\{y_1,\ldots, y_6\}\subset C(e_0,\theta_0)\subset{\bf S}^3,\; $ where $Y$ is an optimal $z$-code. We may assume that
$\theta_1\le\theta_2\le\ldots\le\theta_6.\; $ 
Then from Corollary 3$(i)$ follows that $$\; \theta_0\ge\theta_6\ge\theta_5\ge\arccos{\sqrt{z}}.$$

Let us consider two cases: (a) $\; \theta_0\ge\theta_6\ge\theta_0',\; \, $ 
(b) $\; \theta_0'\ge\theta_6\ge\arccos{\sqrt{z}}.\\ \\$
(a) We have $\; h_6=H(y_0;y_1,\ldots,y_6)=H(y_0;y_1,\ldots,y_5)+\tilde f(\theta_6),$ 
$$H(y_0;y_1,\ldots,y_5)\le h_5=\lambda_5(\psi,\theta_0), \; 
\quad \tilde f(\theta_6)\le\tilde f(\theta_0'). $$
Then $\; h_6 \le \tilde f(\theta_0')+\lambda_5(\psi,\theta_0). \\ \\$
(b) In this case all $\theta_i\le\theta_0',\; $ i.e. $Y\subset C(e_0,\theta_0')$. Since
$$H(y_0;y_1,\ldots,y_5)\le \lambda_5(\psi,\theta_0'), \; 
\quad \tilde f(\theta_6)\le f(-\sqrt{z}), $$
it follows that $\; h_6\le f(-\sqrt{z})+\lambda_5(\psi,\theta_0'). $
\end{proof}

We have proved the following theorem.
\begin{theorem} Suppose $n=4, \;  f\in{\it\Phi}^*(z), \; \sqrt{z}>t_0>z>0$, 
 and $N$ is a positive integer.  Then
$\\ \\ (i)\quad  h_0=f(1),\quad h_1=f(1)+f(-1);\\ \\
(ii) \quad \! h_m=\lambda_m(\psi,\theta_0) 
\le\lambda_m(N,\psi,\theta_0)\; $ for $\; 2\le m\le 5;\\ \\ $  
$(iii)  \; \; h_6\le\max{\{\tilde f(\theta_0')+\lambda_5(\psi,\theta_0),\, f(-\sqrt{z})+\lambda_5(\psi,\theta_0')\}}\; \, \forall \; \, \theta_0'\in [\arccos{\sqrt{z}},\theta_0].$ 
\end{theorem}

\medskip

\noindent{\bf 6-F. Proof of Lemma B.}  
First we show that $f_4\in{\it\Phi}^*(1/2)$ (see Fig. 9). Indeed,
the polynomial $f_4$ has two roots on $[-1,1]$: $t_1=-t_0, \; t_0\approx 0.60794, \; t_2=1/2$; 
$\; f_4(t)\le 0\;$ for $\;t\in [-t_0,1/2],$ and $f_4$ is a monotone decreasing function on the interval $[-1,-t_0].$ The last property holds because there are no zeros of the derivative 
$f'_4(t)$ on  $[-1,-t_0]$. Thus, $f_4\in{\it\Phi}^*(1/2)$.

We have $t_0>0.6058.$ Then Corollary 3$(ii)$ gives $\mu\le 6.$ For calculations of $h_m$ let us apply Theorem 6  with $\psi=\arccos{z}=60^\circ,\; \theta_0=\arccos{t_0}\approx 52.5588^\circ\; $  We get 
$$h_0=f(1)=18.774,\quad h_1=f(1)+f(-1)=24.48.$$ 
$H_2$ achieves its maximum at $\theta_1=30^\circ.$ Then
$$h_2= f(1)+2f(-\cos{30^\circ})\approx 24.8644.$$
For $m=3$ we have
$$h_3=\lambda_3(60^\circ,\theta_0)\approx 24.8345$$ 
at $\theta_3=\theta_0, \; \theta_1=\theta_2\approx30.0715^\circ.$

The polynomial $H_4$ attains its maximum $$\; h_4\approx 24.818\;$$ at the point with $\; \theta_1=\theta_2\approx 30.2310^\circ,\;\; \theta_3=\theta_4\approx 51.6765^\circ,$ and
$$h_5\approx 24.6856$$ 
at $\alpha=60^\circ,\; 
 \theta_1\approx 42.1569^\circ,\; \theta_2=\theta_4=32.3025^\circ, \; \theta_3=\theta_5=\theta_0.$

Let $\theta_0'=50^\circ.$ We have
$\; \tilde f(50^\circ)\approx 0.0906, \; \, \arccos{\sqrt{z}}=45^\circ, \; \, \tilde f(45^\circ)\approx 0.4533,$
$$  \lambda_5(60^\circ,\theta_0)=h_5\approx 24.6856, \quad \; \lambda_5(60^\circ,50^\circ)\approx 23.9181,$$
$$h_6\le \max{\{ \, \tilde f(50^\circ)+h_5,\, \tilde f(45^\circ)+\lambda_5(60^\circ,50^\circ) \, \}}\approx 24.7762<h_2.$$
Thus $\; h_{max}=h_2<25$. Since $(4.2)$, we have $S(X)<25M$.

\begin{center}
\begin{picture}(320,200)(-160,-110)
\thinlines
\put(-135,-80){\line(0,1){160}}
\put(135,-80){\line(0,1){160}}
\put(-135,-80){\line(1,0){270}}
\put(-135,80){\line(1,0){270}}
\put(-135,-60){\line(1,0){270}}

\thicklines
\qbezier (-135,54)(-132,51)(-129,45)
\qbezier (-129,45)(-126,37)(-123,27)
\qbezier (-123,27)(-120,18)(-117,8)
\qbezier (-117,8)(-114,-2)(-111,-11)
\qbezier (-111,-11)(-108,-19)(-105,-27)
\qbezier (-105,-27)(-99,-40) (-93,-49)
\qbezier (-93,-49)(-90,-52)(-87,-55)
\qbezier (-87,-55)(-84,-57)(-81,-59)
\qbezier (-81,-59)(-78,-60)(-75,-61)
\qbezier (-75,-61)(-69,-62)(-53,-62)
\qbezier (-53,-62)(-45,-61)(-37,-61)
\qbezier (-37,-61)(-10,-61)(15,-61) 
\qbezier (15,-61)(23,-62)(30,-63)
\qbezier (30,-63)(38,-65)(45,-67)
\qbezier (45,-67)(53,-69)(60,-71) 
\qbezier (60,-71)(68,-72)(75,-71) 
\qbezier (75,-71)(83,-68)(90,-60)
\qbezier (90,-60)(98,-50)(105,-35)
\qbezier (105,-35)(113,-16)(120,8)
\qbezier (120,8)(128,36)(134,69)

\thinlines
\multiput (-120,-80)(15,0){17}%
{\line(0,1){2}}
\multiput (-135,-40)(0,20){6}%
{\line(1,0){2}}
\put(-143,-90){$-1$}
\put(-119,-90){$-0.8$}
\put(-89,-90){$-0.6$}
\put(-59,-90){$-0.4$}
\put(-29,-90){$-0.2$}
\put(13,-90){$0$}
\put(40,-90){$0.2$}
\put(70,-90){$0.4$}
\put(100,-90){$0.6$}
\put(130,-90){$0.8$}
\put(-143,-62){$0$}
\put(-143,-42){$1$}
\put(-143,-22){$2$}
\put(-143,-2){$3$}
\put(-143,18){$4$}
\put(-143,38){$5$}
\put(-143,58){$6$}
\put(-150,-82){$-1$}
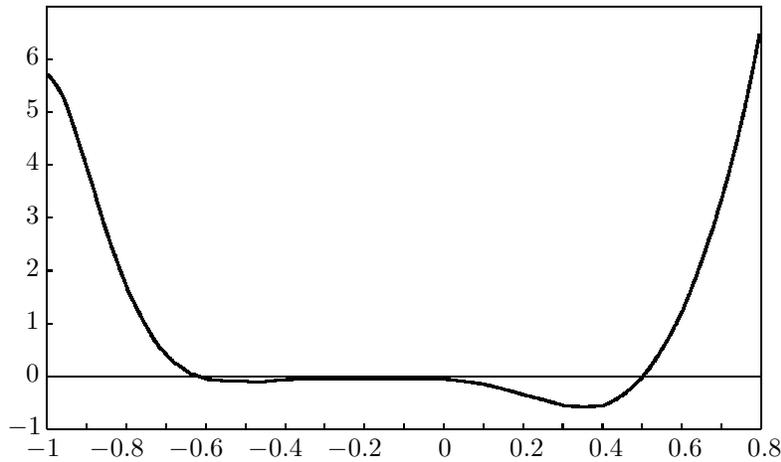
\put(-78,-110){Fig. 9. The graph of the function $f_4(t)$}

\end{picture}
\end{center}

\section {Concluding remarks }
This extension of the Delsarte method can be applied to other dimensions and spherical $\psi$-codes. 
 The most interesting application is a new proof for the Newton-Gregory problem, $k(3)<13.$ In dimension three computations of $h_m$ are technically  much more easier than for $n=4$ (see \cite{Mus2}).

Let 
$$f(t) = \frac{2431}{80}t^9 - \frac{1287}{20}t^7 + \frac{18333}{400}t^5 + \frac{343}{40}t^4 - \frac{83}{10}t^3 - \frac{213}{100}t^2+\frac{t}{10} - \frac{1}{200}. $$
Then $f\in{\it\Phi}^*(1/2),\; t_0\approx 0.5907, \; \mu(3,1/2,f)=4,$ and $h_{max}=h_1=12.88.$
The expansion of $f$ in terms of Legendre polynomials $P_k=G_k^{(3)}$ is
$$f = P_0 + 1.6P_1 + 3.48P_2 + 1.65P_3 + 1.96P_4 + 0.1P_5 + 0.32P_9.$$
Since $c_0=1,\; c_i\ge 0,$ we have $  k(3)\le h_{max}=12.88<13.$

Direct application of the method developed in this paper, presumably could lead to some improvements in the upper bounds on kissing numbers in dimensions 9, 10, 16, 17, 18 given in \cite[Table 1.5]{CS}. (``Presumably" because the equality $\; h_{max}=E\; $ is not proven yet.)

In 9 and 10 dimensions Table 1.5 gives: $\\306\le k(9)\le 380,\quad 500\le k(10)\le 595.$\\
Our method gives:\\
$n=\;\,9:\; \deg{f}=11,\; E=h_1=366.7822,\; t_0=0.54;$\\
$n=10:\; \deg{f}=11,\; E=h_1=570.5240,\; t_0=0.586$.\\
For these dimensions there is a good chance to prove that $\\k(9)\le 366,\; k(10)\le 570.$

From the equality $k(3)=12$ follows that $\varphi_3(13)<60^\circ.$ 
The method gives \\ $\varphi_3(13)<59.4^\circ$ ($\deg{f}=11$).
The lower bound on $\varphi_3(13)$ is $57.1367^\circ $ \cite{FeT}. Therefore, we have $57.1367^\circ\le\varphi_3(13)<59.4^\circ.$

Using our approach it can be proven that $\varphi_4(25)<59.81^\circ, \; \varphi_4(24) < 60.5^\circ.$ 
That improve the bounds:
$$\varphi_4(25)<60.79^\circ, \;\; \varphi_4(24) < 61.65^\circ \; \cite{Lev2} \; (\mbox{cf. } \cite{Boyv});
\;\; \varphi_4(24) < 61.47^\circ \; \cite{Boyv};$$
$$ \varphi_4(25)<60.5^\circ, \quad \varphi_4(24) < 61.41^\circ \; \cite{AB2}.$$
Now in these cases we have  
$$\quad 57.4988^\circ\ < \varphi_4(25) < 59.81^\circ,
\quad 60^\circ \le \varphi_4(24) <  60.5^\circ. \footnote{The long-standinding conjecture: The maximal kissing arrangment in four dimensions is unique up to isometry (in other words, is the ``24-cell"), and
$\varphi_4(24)=60^\circ$.}
$$

\medskip

However, for $n=5,6,7$ direct use of this extension of the Delsarte method doesn't give better upper bounds on $k(n)$ than Odlyzko-Sloane's bounds \cite{OdS}. It is an interesting problem to find better methods.

\pagebreak

\centerline{\bf\Large Appendix. An algorithm }
\centerline{\bf\Large for computation suitable polynomials $f(t)$}

\medskip

\medskip

In this Appendix is presented an algorithm for computation  ``optimal" \footnote{Open problem: is it true that for given $t_0, d$ this algorithm defines $f$ with minimal $h_{max}$?}  polynomials $f$ such that 
$f(t)$ is a monotone decreasing function on the interval $[-1,-t_0],$ and  
$f(t)\le 0 \; \mbox{ for } \; t\in [-t_0,z], \quad t_0>z\ge 0$. This algorithm based on our knowledge about optimal arrangement of points $y_i$ for given $m$. Coefficients $c_k$ can be found via discretization and linear programming; such method had been employed already by Odlyzko and Sloane \cite{OdS} for the same purpose.

Let us have a polynomial $f$ represented in the form $f(t)=1+\sum\limits_{k=1}\limits^d c_kG_k^{(n)}(t)$. We have the following constraints for $f$: (C1) $\;\; c_k\ge 0,\;\; 1\le k\le d$;\\ (C2) $\; f(a)>f(b)\;$ for $\; -1\le a<b\le -t_0$;\quad (C3) $\; f(t) \le 0\;$ for $\; -t_0\le t\le z.$

We do not know $e_0$ where $H_m$ attains its maximum, so for evaluation of $h_m$ let us use $e_0=y_c,$ where $y_c$ is the center of $\Delta_m.$ All vertices $y_k$ of $\Delta_m$ are at the distance of $\rho_m$ from $y_c,$ where $$\cos{\rho_m}=\sqrt{(1+(m-1)z)/m}.$$

When $m=2n-2, \; \Delta_m$ presumably is  a regular $(n-1)$-dimensional cross-polytope.\footnote {It is also an open problem.} In this case $\; \cos{\rho_m}=\sqrt{z}.$ 

Let $I_n=\{1,\ldots,n\}\bigcup \{2n-2\}, \;\; m\in I_n, \;\; b_m=-\cos{\rho_m},\; $ then \\ $H_m(y_c)=f(1)+mf(b_m).\;\; $ If $F_0$ is such that $H(y_0;Y) \le E=F_0+f(1),\;$ then (C4) $\; f(b_m)\le F_0/m, \;\; m\in I_n.\;$ Note that $E=F_0+1+c_1+\ldots+c_d=F_0+f(1)$ is a lower estimate of $h_{max}$. 
A polynomial $f$ that satisfies (C1-C4) and gives the minimal $E$ 
can be found by the following 

\medskip

\centerline{\bf Algorithm.} 

\medskip

\noindent Input: $\; n,\; z,\; t_0,\; d,\; N.$

\noindent Output: $\; c_1,\ldots, c_d,\; F_0, \; E.$

\medskip

\noindent {\it First} replace (C2) and (C3) by a finite set of inequalities at the points\\ $a_j=-1+\epsilon j,\;\; 0\le j \le N, \;\; \epsilon=(1+z)/N:$ 

\medskip

\noindent {\it Second} use linear programming to find $F_0, c_1,\ldots, c_d$ so as to minimize \\ 
$E-1=F_0+\sum\limits_{k=1}\limits^dc_k\;\;$ subject to the constraints
$$c_k\ge 0,\quad 1\le k\le d;\qquad \sum\limits_{k=1}\limits^dc_kG_k^{(n)}(a_j)\ge \sum\limits_{k=1}\limits^dc_kG_k^{(n)}(a_{j+1}), \quad a_j\in [-1,-t_0];$$
 $$1+\sum\limits_{k=1}\limits^dc_kG_k^{(n)}(a_j)\le 0,\quad a_j \in [-t_0,z];\quad 
1+\sum\limits_{k=1}\limits^dc_kG_k^{(n)}(b_m)\le F_0/m,\quad m\in I_n.$$ 

Let us note again that $E \le h_{max}$, and $E = h_{max}$ only if $h_{max} = H_{m_0}(y_c)$ for some  $m_0\in I_n.$



\begin{thebibliography}{99}

\bibitem{Ans}
K. Anstreicher, The thirteen spheres: A new proof, Discrete and Computational Geometry, {\bf 31}(2004), 613-625.
\bibitem{AB1}
V.V. Arestov and A.G. Babenko, On Delsarte scheme of estimating the contact numbers, 
Proc. of the Steklov Inst. of Math. {\bf 219} (1997), 36-65. 

\bibitem{AB2}
V.V. Arestov and A.G. Babenko, Estimates for the maximal value of the angular code distance for 24 and 25 points on the unit sphere in ${\bf R}^4$, Math. Notes, {\bf 68} (2000), 419-435.

\bibitem{Boyv}
P.G. Boyvalenkov, D.P. Danev and S.P. Bumova, Upper bounds on the minimum distance of spherical codes, IEEE Trans. Inform. Theory, {\bf 42}(5), 1996, 1576-1581. 
\bibitem{Bor1}
K. B\"or\"oczky, Packing of spheres in spaces of constant curvature, Acta Math. Acad. Sci. Hung.  {\bf 32} (1978), 243-261.
\bibitem{Bor}
K. B\"or\"oczky, The Newton-Gregory problem revisited, Proc. Discrete Geometry, Marcel Dekker, 2003,  103-110.
\bibitem{Car}
B.C. Carlson, Special functions of applied mathematics, Academic Press, 1977.
\bibitem{Cas}
B. Casselman, The difficulties of kissing in three dimensions, Notices Amer. Math. Soc., {\bf 51}(2004), 884-885. 
\bibitem{CS}
J.H. Conway and N.J.A. Sloane, Sphere Packings, Lattices, and Groups, New York, Springer-Verlag, 1999 (Third Edition).
\bibitem{Cox}
H.S.M. Coxeter, An upper bound for the number of equal nonoverlapping spheres that can touch another of the same size, Proc. of Symp. in Pure Math. AMS, {\bf 7} (1963), 53-71 = Chap. 9 of
H.S.M. Coxeter, Twelve Geometric Essays, Southern Illinois Press, Carbondale Il, 1968.
\bibitem{Dan}
L. Danzer, Finite point-sets on ${\bf S}^2$ with minimum distance as large as possible, Discr. Math., {\bf 60} (1986), 3-66. 
\bibitem{DGK}
L. Danzer, B. Gr\"unbaum, and V. Klee. Helly's theorem and its relatives. Proc. Sympos. 
Pure Math., vol. 7, AMS, Providence, RI, 1963, pp. 101-180.
\bibitem{Del1}
Ph. Delsarte, Bounds for unrestricted codes by linear programming, Philips Res. Rep., {\bf 27}, 1972, 272-289.
\bibitem{Del2}
Ph. Delsarte, J.M. Goethals and J.J. Seidel, Spherical codes and designs, Geom. Dedic., {\bf 6}, 1977, 363-388.
\bibitem{Erd}
A. Erd\'elyi, editor, Higher Transcendental Function, McGraw-Hill, NY, 3 vols, 1953, Vol. II,
Chap. XI.
\bibitem{FeT} 
L. Fejes T\'oth, Lagerungen in der Ebene, auf der Kugel und in Raum, Springer-Verlag, 1953; Russian translation, Moscow, 1958.

\bibitem{Hales}
T. Hales, The status of the Kepler conjecture, Mathematical Intelligencer {\bf 16}(1994), 47-58.

\bibitem{Hop}
R. Hoppe, Bemerkung der Redaction, Archiv Math. Physik (Grunet) {\bf 56} (1874), 307-312. 
\bibitem{Hs1}   
W.-Y. Hsiang, The geometry of spheres, in Differential Geometry (Shanghai,1991), Word Scientific, River Edge, NJ, 1993, pp. 92-107. 
\bibitem{Hs}
W.-Y. Hsiang, Least Action Principle of Crystal Formation of Dense Packing Type and Kepler's Conjecture, World Scientific, 2001. 
\bibitem{Kab}
G.A. Kabatiansky and V.I. Levenshtein, Bounds for packings on a sphere and in space,
Problems of Information Transmission, {\bf 14}(1), 1978, 1-17. 
\bibitem{Lee}
J. Leech, The problem of the thirteen spheres, Math. Gazette {\bf 41} (1956), 22-23. 
\bibitem{Lev2}
V.I. Levenshtein, On bounds for packing in $n$-dimensional Euclidean space, Sov. Math. Dokl. 
{\bf 20}(2), 1979, 417-421.


\bibitem{Ma}
H. Maehara, Isoperimetric theorem for spherical polygons and the problem of 13 spheres, Ryukyu Math. J.,
{\bf 14} (2001), 41-57.

\bibitem{Mus}
O.R. Musin, The problem of the twenty-five spheres, Russian Math. Surveys, {\bf 58}(2003), 794-795. 

\bibitem{Mus2}
O.R. Musin, The kissing problem in three dimensions,  Discrete Comput. Geom., {\bf 35} (2006), 375-384.

\bibitem{OdS}
A.M. Odlyzko and N.J.A. Sloane, New bounds on the number of unit spheres that
can touch a unit sphere in $n$ dimensions, J. of Combinatorial Theory
A26(1979), 210-214.

\bibitem{PZ}
F. Pfender and G.M. Ziegler, Kissing numbers, sphere packings, and some unexpected proofs, Notices Amer. Math. Soc., {\bf 51}(2004), 873-883. 

\bibitem{Scho}
I.J. Schoenberg,  Positive definite functions on spheres, Duke Math. J., 
{\bf 9} (1942), 96-107.
\bibitem{SvdW1}
K. Sch\"utte and B.L. v. d. Waerden, Auf welcher Kugel haben 5,6,7,8 oder 9
Punkte mit Mindestabstand 1 Platz? Math. Ann. {\bf 123} (1951), 96-124.
\bibitem{SvdW2}
K. Sch\"utte and B.L. van der Waerden, Das Problem der dreizehn Kugeln, Math. Ann. {\bf 125} (1953), 325-334.
\bibitem{Wyn}
A.D. Wyner, Capabilities of bounded discrepancy decoding, Bell Sys. Tech. J. {\bf 44} (1965), 1061-1122. 


\end{thebibliography}
\end{document}